\def\WHO{nobody} 
\def\version{25.10.2024}
\def\users{final-layout}   
\documentclass[12pt]{article} 
\usepackage{upgreek}
\usepackage{xcolor}
\usepackage{bm,amsmath,amsthm,hyperref,amsfonts,amssymb,color}
\usepackage{mathrsfs} 
\usepackage{cite}

\usepackage{ifthen}
\ifthenelse{\equal{\users}{final-layout}}{}{
\usepackage{fancyhdr}
\pagestyle{fancy}
\headheight=28pt\headwidth=17cm
\definecolor{gray}{gray}{0.5}
\rhead{\color{gray}gravitational differentiation of planets\\
A. Mielke \& T.Roub\'\i\v cek \& U.Stefanelli}
\chead{}
\lhead{Version \version, file: \jobname.tex,
\\
compiled:\number\day.\number\month.\number\year\ at
\the\hour:\ifnum\minute<10 0\fi\the\minute\ h\ 
\underline{\Large\color{red}\WHO's working on!}
}
}
 
\newcount\hour \newcount\minute
\hour=\time
\divide \hour by 60
\minute=\time
\loop \ifnum \minute > 59 \advance \minute by -60 \repeat
  
\definecolor{labelkey}{rgb}{1.,.2,0.} 
 
\usepackage[normalem]{ulem}
\usepackage{ifthen}
\usepackage{color} 

\usepackage{graphicx}
\usepackage{psfrag} 
\newcounter{myfigure}
\newenvironment{my-picture}[3]{\refstepcounter{myfigure}\label{#3}\setlength{\unitlength}{1cm}\begin{picture}(#1,#2)}{\end{picture}}

\ifthenelse{\equal{\users}{final-layout}}{

	\newcommand{\COMMENT}[1]{}
	\newcommand{\COMMENTGT}[1]{}
	\newcommand{\TODO}[1]{}
	\newcommand{\INTERNAL}[1]{}
	\newcommand{\QUESTION}[1]{}
	\newcommand{\DELETE}[1]{}

	\newcommand{\REM}[1]{\marginpar{\bfseries\tiny{\color{blue}}}}
    \newcommand{\MARGINOTE}[1]{}

}
{
	
	\newcommand{\COMMENT}[1]{{\color{red}\uuline{#1}\color{black}}}
	\newcommand{\COMMENTGT}[1]{{\hfill\large\color{red}***{#1}***\color{black}\hfill}\\}
	\newcommand{\TODO}[1]{{\color{red}\uuline{#1}\color{black}}}
	\newcommand{\INTERNAL}[1]{\footnote{#1}}
	\newcommand{\QUESTION}[1]{{\color{brown}\uuline{#1}\color{black}}}
	\newcommand{\DELETE}[1]{{\color{red}\sout{#1}\color{black}}}

	\newcommand{\REM}[1]{\marginpar{\bfseries\tiny{\color{blue}#1}}}
\newcommand{\MARGINOTE}[1]{\marginpar{\color{red}\tiny\texttt{#1}}}

}

\newcommand\DT[1]{\mathchoice
                 {{\buildrel{\hspace*{.1em}\text{\LARGE.}}\over{#1}}}
                 {{\buildrel{\hspace*{.1em}\text{\LARGE.}}\over{#1}}}
                 {{\buildrel{\hspace*{.1em}\text{\Large.}}\over{#1}}}
                 {{\buildrel{\hspace*{.1em}\text{\large.}}\over{#1}}}}
\newcommand\pdt[1]{\frac{\partial{#1}}{\partial t}} 
\newcommand{\lineunder}[2]{\LU{\begin{array}[t]{c}\underbrace{#1}\vspace*{.5em}\end{array}}{\mbox{\footnotesize\rm #2}}}
\newcommand{\LU}[2]{\begin{array}[t]{c}#1\vspace*{-1em}\\_{#2}\end{array}}
\newcommand{\linesunder}[3]{\LSU{\begin{array}[t]{c}\underbrace{#1}\vspace*{.5em}\end{array}}{\mbox{\footnotesize\rm #2}}{\mbox{\footnotesize\rm #3}}}
\newcommand{\LSU}[3]{\begin{array}[t]{c}#1\vspace*{-1em}\\_{#2}\vspace*{-.5em}\\_{#3}\end{array}}
\newcommand{\morelinesunder}[4]{\LSUU{\begin{array}[t]{c}\underbrace{#1}\vspace*{.5em}\end{array}}{\mbox{\footnotesize\rm #2}}{\mbox{\footnotesize\rm #3}}{\mbox{\footnotesize\rm #4}}}
\newcommand{\LSUU}[4]{\begin{array}[t]{c}#1\vspace*{-1em}\\_{#2}\vspace*{-.5em}\\_{#3}\vspace*{-.5em}\\_{#4}\end{array}}
\newcommand{\Item}[2]{\parbox[t]{.055\textwidth}{#1}\hfill%
      \parbox[t]{.945\textwidth}{#2}\vspace*{.8mm}} 
\newcommand{\divS}{\mathrm{div}_{\scriptscriptstyle\textrm{\hspace*{-.1em}S}}^{}}
\newcommand{\nablaS}{\nabla_{\scriptscriptstyle\textrm{\hspace*{-.3em}S}}^{}}

\def\Vdots{\!\mbox{\setlength{\unitlength}{1em}
\begin{picture}(0,0)
\put(-.07,0){.}
\put(-.07,.3){.}
\put(-.07,.6){.}
\end{picture}\hspace*{.2em}}}
\usepackage{mathrsfs}   
\usepackage{eucal}      

  \def\bbI{{\bm{I}}}
 \def\bbK{{\mathbb K}} 
\def\bbM{{\mathbb M}}  
\def\bbP{{\mathbb P}} \def\bbQ{{\mathbb Q}}

\def\FG{\boldsymbol}
  \def\cc{{\FG c}} 
 \def\ee{{\FG e}}  
\def\gg{{\FG g}}  
\def\jj{{\FG j}}   
 \def\nn{{\FG n}}  
  \def\rr{{\FG r}} 

\def\vv{{\FG v}}  \def\xx{{\FG x}} 
\def\yy{{\FG y}}  
   
\def\DD{{\FG D}} 
\def\FF{{\FG F}}

 \def\NN{{\FG N}}  
   
 \def\TT{{\FG T}} 
  \def\XX{{\FG X}}

\newcommand{\R}{\mathbb R}
 
\newcommand{\Nabla}{{\nabla}}

\newcommand{\EE}{{\bm e}}

\newcommand{\pl}{\partial}
\newcommand{\eq}[1]{(\ref{#1})}

\renewcommand{\d}{\mathrm d}  
\newcommand{\barOmega}{\hspace*{.2em}{\overline{\hspace*{-.2em}\varOmega}}}
\newcommand{\ubarOmega}{\hspace*{.05em}{\overline{\hspace*{-.05em}\Omega\hspace*{-.15em}}}\hspace*{.15em}}

 \newcommand{\bulet}{\text{\footnotesize$\,\bullet\,$}}

\newcommand{\rhoR}{\varrho_\text{\sc r}}
\newcommand{\rhoRxi}{\varrho_\text{\sc r}^{\bm\xi}}

\newcommand{\rhoREF}{\varrho_\text{\rm ref}}

\newcommand{\W}{w}
\newcommand{\COUPLING}{\gamma}
\newcommand{\wh}{\widehat}
\newcommand{\wt}{\widetilde}

\newtheorem{theorem}{Theorem}[section]

\newtheorem{definition}[theorem]{Definition}

\newtheorem{remark}[theorem]{Remark}

\numberwithin{equation}{section}
 
\newcommand{\AAA}{\color{black}}
\newcommand{\UUU}{\color{black}}

\newcommand{\TTT}{\color{black}}
\newcommand{\SSS}{\color{black}}     
\newcommand{\FFF}{\color{black}} 
\newcommand{\EEE}{\color{black}}

\def\Universe{\R^3}

\begin{document}

\allowdisplaybreaks

\begin{center}

\noindent{\LARGE\bf 
  A model of gravitational differentiation of compressible
  self-gravitating planets}

\bigskip\bigskip\bigskip

\noindent{\bf Alexander Mielke}\footnote{Weierstrass Institute for Applied
Analysis and Stochastics, Mohrenstra\ss{}e 39, 10117 Berlin,
Germany.\\ \hspace*{1.6em} Email: \texttt{alexander.mielke@wias-berlin.de}}%
\footnote{Institut  f\"ur  Mathematik,  Humboldt--Universit\"at zu
  Berlin,  Rudower Chaussee 25,  D-12489  Berlin, Germany.},
         {\bf Tom\'a\v s Roub\'\i\v cek}\footnote{Mathematical Institute, Charles University,
Sokolovsk\'a 83, CZ-186~75~Praha~8,  Czech Republic.\\\hspace*{1.6em} Email: \texttt{tomas.roubicek@mff.cuni.cz}}
\footnote{Institute of Thermomechanics, Czech Academy of Sciences,
Dolej\v skova 5, CZ-18200~Praha~8, Czech Republic},
and {\bf Ulisse Stefanelli}\footnote{Faculty of Mathematics, University of
  Vienna, Oskar-Morgenstern-Platz 1, 1090 Vienna, Austria and Vienna Research Platform on Accelerating Photoreaction Discovery, University of Vienna, W\"ahringerstrasse 17, A-1090 Vienna, Austria.\\\hspace*{1.6em} Email: \texttt{ulisse.stefanelli@univie.ac.at}}
\footnote{Istituto di Matematica Applicata e Tecnologie Informatiche
  {\it E. Magenes} - CNR, v. Ferrata 1, I-27100 Pavia, Italy.}
  
\bigskip\bigskip\bigskip\bigskip

{\small\bf Abstract}.

\vspace*{.3em}

\begin{minipage}[t]{36em}{\small
\baselineskip=13pt

 We present a  dynamic model  for  inhomogeneous 
viscoelastic media  at   finite strains.  The model features a
Kelvin-Voigt rheology, and  includes  a  self-generated gravitational field in
the  actual evolving configuration.  In particular, a fully Eulerian approach
is adopted. We specialize  the model to viscoelastic (barotropic) fluids
and prove   existence and a certain regularity of global weak solutions
by a Faedo-Galerkin semi-discretization  technique. Then, an extension to
multi-component chemically reacting viscoelastic fluids based on a
phenomenological approach  by  Eckart and Prigogine, is advanced and studied.
The model is inspired by planetary geophysics. In particular, it describes
gravitational differentiation of inhomogeneous planets and moons, possibly
undergoing volumetric phase transitions.

\medskip

\medskip


}
\end{minipage}
\end{center}

\bigskip\bigskip

\baselineskip=16pt

\def\GRAVPOT{V}
\def\GRAVCONST{G}

\section{Introduction}

Self-gravitating  inhomogeneous
media  provide a rich class of interesting problems  
in continuum mechanics,  a  prominent application  being  planetary
geophysics.  A detailed understanding of the processes of
planetary formation and early evolution is currently available. In
addition, the dynamics of the interiors of planets and their moons is
relatively well understood, in particular in relation with the Solar
planetary  system and the planet Earth, cf.\
\cite{BadWal15EEAD,Cond16EEPS,Gery19INGM,GerYue07RCMM,ShTuOl04MCEP}.

After a relatively short  time  (tens or hundreds of millions of years)
of accretion from a rather homogeneous stellar disc of dust, meteoroids, and
asteroids, on  
a much longer period (billions of years) planets {\it differentiate}.  
 Self-gravitation drives the dynamics inside the
mantle, eventually leading to the formation of a {\it core-mantle}
structure. This  occurs as effect of the different
densities of the mantle constituents, i.e., the heavier media
(metals), which  are strongly attracted towards the planet core, and
the lighter ones
(silicates and volatile elements forming liquid oceans or gaseous
atmospheres),  which are subjected to buoyancy. The onset of so-called
Rayleigh-Taylor instability of interfaces between media with  
different densities can be observed. A similar evolution happens also
in the development of big moons of planets.

Self-gravitation also governs the subsequent evolution of planets and
moons interiors. In particular, gravity is responsible for the
formation of plumes and slabs in the mantle, which are accompanied by
various volumetric {\it phase transitions} and related buoyancy effects.
Silicates in the Earth mantle undergo several quite sharp volumetric
transitions. This happens at  
pressures of about 14 GPa (olivine transforms to wadsleyite) and
23 GPa (spinel transforms to perovskite and magnesiow\"ustite), which in
Earth's mantle  occur at a depth of 410\;km and 660\,km, respectively, cf.,
e.g., \cite{Chri95EPTM,Cond16EEPS,Gery19INGM,HelWoo01EM,TSGS94EMPT}. These
volumetric transformations \TTT may account for hysteresis and \EEE are related
with the loss of  strict 
convexity \TTT of the stored
energy $\phi_{\rm ref}$ \EEE as a function of the determinant  $J=\det\FF$ of the
deformation gradient, i.e., {\it  loss of  polyconvexity} in terms of
$\FF$ itself.

 The aim of this paper is to present and analyze a dynamic model
able to capture the basic features of planets and moons
differentiation under self-gravitation. In particular, we model an inhomogeneous viscoelastic
medium at finite strains by assuming the so-called Kelvin-Voigt
rheology and including the effects of the self-generated gravitational
field. A crucial aspect of our approach is that it is 
formulated in the actual evolving configuration, making it fully
Eulerian. 
 Mainly  for analytical reasons  (but see 
Remarks~\ref{rem-anisotherm0} and \ref{rem-anisotherm1}  for a
discussion),  we restrict ourselves to the discussion of the
isothermal case. Note nonetheless that  
thermal effects and temperature dependence of material parameters  play a
vital role in the  evolution of planets' and moons'
interiors.  We refrain however to discuss thermal effects here,
as these would require a much longer 
tractate.   
We specifically study the case of barotropic fluids,
as well as the multicomponent case. Both in the single component and
in the multicomponent case, we are able to prove the existence of weak
solutions by means of a Galerkin approximation technique
\FFF on long time scales, cf.\ Remark~\ref{rem-long-time} below\EEE.

Before  starting our discussion of the evolutionary model, 
 let us mention some available contributions on the
equilibrium of self-gravitating systems. By assuming spherical
symmetry, a collection of results, together with an
historical account, are provided in the monograph 
\cite{MulWei16SDES}. 
Existence of equilibria of a hyperelastic solid under self-gravitation
has been investigated in \cite{CalLeo15GSSG}. Some alternative
Eulerian analysis in case of spheres and multiple and stratified spherical shells is in \cite{Alho19,Alho20}.
Again in the spherical setting, the phenomenon of gravitational
collapse and the possible existence of multiple equilibria are
discussed in \cite{Jia19}, together with the ensuing bifurcation
dynamics.

The case of nonlinear elastic models for polytropic fluids are
considered in \cite{Cal22} both in the
Lagrangian and the Eulerian setting. In the  case of the  stationary
Navier-Stokes under barotropic pressure  with  $p(\rho) \approx 
a\rho^\gamma$ for $\gamma>4/3$,  existence of equilibria  has been
established in
\cite{sec94}, and the stability of radially symmetric solutions and
their free boundary have been studied in \cite{StrZaj99LESF}.
In addition, a number of contributions explore the relativistic
setting, again in the spherical case, see \cite{Alho17,And14} among many
others. 

In the dynamic case, the reference setting is that of Navier-Stokes-Poisson
systems. Note that these are parted into two distinct classes, depending on the
repulsive versus attractive nature of the Poisson subsystem. These indeed
correspond to the  modeling of  electrically  charged versus
self-gravitating fluids.  The literature in the repulsive case is extensive,
see, for instance \cite{CDGS21DMVS, Dona03LGEC, Dona08, He20,
  LiMaZh10ODRC,Tan10,Li10,Zhe12}.   Concerning  the attractive case,
which  is the case of our interest,  we start by mentioning the local
existence result in \cite{StrZaj99LESF}, as well as the global existence
result for weak solutions on bounded domains  with pressure law 
$p(\rho)=a\rho^\gamma$, $\gamma>3/2$, in \cite{Koba08}. Existence in an
external domain is tackled in \cite{DucFei04DGS,DFPS01EGFB,JiaTan09CBTA} for
$p(\rho)=a \rho^\gamma$, $\gamma>3/2$, and in \cite{Duco04} for $p(\cdot)$
non-monotone.  Global existence in $\mathbb{R}^n$ for $4/3 < \gamma \leq 3/2$
and for radially symmetric initial data has been proved in \cite{Jia10}. The
stability and stabilization of spherical solutions in an external domain is
treated in \cite{Duco05}. Eventually, a result on strong-weak uniqueness is
in \cite{Basa22WSUP}, and the existence of an absorbing set for large times
has been obtained in \cite{JiaTan09CBTA} for $\gamma> 5/3$ and later in
\cite{Guo12} for $\gamma>3/2$. Local existence for self-gravitating inviscid
liquid bodies with varying shape was shown in \cite{StrZaj99LESF} and,
in the degenerate (inviscid)  case of the Euler-Poisson system, 
in \cite{GiLiLu20LWPM}.  

The novelty of our paper relies on the treatment of general 
{\it inhomogeneous}, as well {\it multicomponent} materials. Note that the
discussion of inhomogeneous or multicomponent materials is instrumental to
the description of differentiation dynamics in planets' mantle. Our way of
modeling  the inhomogeneity of the material  relies on the use of the  
{\it reference mapping} $\bm\xi(t,\cdot)$ mapping the Eulerian or spatial
point $\bm x\in \bm y(t,\Omega)$ back into the reference domain
configuration $\Omega$. By studying the transport of $\bm\xi$ along with
the velocity field $\bm v(t,\bm x)$ we can trace the inhomogeneities that
are imprinted into the material at the reference point $\bm\xi(t,\bm x)$.
Using a suitable hyperviscosity (cf., \ref{T+D}), we are able to derive
the necessary regularity properties for $\bm\xi$ to treat quite general
material laws,  see \eqref{ass}. In particular, our analysis covers the
case of a general non-monotone barotropic pressure with growth
$\sim\rho^{\gamma}$ with $\gamma>
\TTT2\EEE$; cf., assumption \eqref{ass-phi} with
Remark~\ref{rem-state-eq} below. 

The plan of the paper is as follows.  Section~\ref{sec-kinem} recalls basic
notations and  standard concepts from kinematic of finite-strain continuum
mechanics. In Section~\ref{sec-single},  we introduce the model for
self-gravitating inhomogeneous viscoelastic bodies, and its energetics is
discussed. Section~\ref{sec-fluids} then specifies the model to
visco-elastic Navier-Stokes fluids with Kelvin-Voigt rheology in the
volumetric part. In Section~\ref{sec-anal}, we  present a result on
existence and regularity of weak  solutions to an initial-boundary-value
problem for the inhomogeneous self-gravitating system. In order to achieve  
this goal, we  perform an approximation and Galerkin discretization   
of (most of) the equations. Eventually, in Section~\ref{sec-multi}, we
present and analyze a {\it multi-component} version of the previous model, 
employing the Eckart-Prigogine approach.

\section{Kinematics at finite strains}\label{sec-kinem}

\def\RRR{\text{\sc r}}

Let us start by recalling some basic notion from the general
theory of large (or finite) deformations in continuum mechanics,
cf., e.g., \cite{GuFrAn10MTC,Mart19PCM}. 

 The basic quantity describing the time-dependent evolution of a
deformable body is the (referential)  
{\it deformation}  or {\it motion} 
$\yy:I{\times} \Omega  \to\R^3$, 
where $I=[0,T]$ and $T>0$ is some  given  final time. For all time instants
$t\in I$, the deformation $\yy(t, \cdot )$
maps the {\it reference configuration}
$ \Omega  \subset\R^3$ of the deformable body to its {\it actual}
configuration
$\yy(t, \Omega)$, a subset of the {\it physical space} $\R^3$. In what
follows, we indicate {\it referential coordinates} by $\XX\in
\Omega$ 
and {\it actual coordinates} by $\xx\in\R^3$. By assuming 
$\yy(t,\cdot)$ to be globally invertible, we indicate  its  inverse by
$\bm{\xi}(t,\cdot)=\yy^{-1}(t,\cdot):\yy(t,  \Omega  )\to 
\Omega  $.  Such $\bm\xi$ is usually referred to as  {\it
  return} or {\it reference mapping}, or sometimes {\it inverse motion}. 
 Two  basic kinematic quantities  
are the Lagrangian velocity $\vv_\RRR=\pdt{}\yy$ and the Lagrangian
deformation gradient $\FF_\RRR=\Nabla_{\!\XX}^{}\yy$.  Starting
from these, one defines  
the Eulerian velocity $\vv(t,\xx)=\vv_\RRR(t,\bm\xi(t,\xx))$ and the Eulerian
deformation gradient $\FF(t,\xx)=\FF_\RRR(t,\bm\xi(t,\xx))$.
Here and  throughout  the  
article, having the Eulerian velocity at disposal,
we use the dot-notation
$(\cdot)\!\DT{^{}}=\pdt{}+\vv{\cdot}\nabla_\xx$  to indicate  the {\it convective time
derivative} applied to scalars or, componentwise, to vectors or
tensors. 
 The velocity gradient $\Nabla\vv$ fulfills  $\Nabla\vv=\nabla_{\!\XX}^{}\vv\nabla_{\!\xx}^{}\XX=\DT\FF\FF^{-1}$,
where we used the chain-rule 
and  the fact that 
$\FF^{-1}=(\nabla_{\!\XX}^{}\xx)^{-1}=\nabla_{\!\xx}^{}\XX$. 
This gives the {\it  evolution-and-transport}  equation  for the
deformation gradient 
\begin{align}
\DT\FF=(\nabla\vv)\FF\,.
  \label{ultimate}\end{align} 
From this, we also obtain the evolution-and-transport equation for Jacobian
$\det\FF$ and its  reciprocal $1/\det\FF$, namely,   
\begin{align}\nonumber\\[-2.5em]
\DT{\overline{\det\FF}}=(\det\FF){\rm div}\,\vv\ \ \ \ \text{ and }\ \ \ \
\DT{\overline{\!\!\!\bigg(\frac1{\det\FF}\bigg)\!\!\!}}\ =-\frac{{\rm div}\,\vv}{\det\FF}\,.
\label{DT-det}\end{align} 
The return mapping $\bm\xi$ satisfies the transport equation
\begin{align} 
\DT{\bm\xi}=\bm0\, 
\label{transport-xi}
\end{align}
 which simply expresses the fact that the material properties encoded in
the material point $\bm X=\bm\xi(t,\bm x)$ move along with the particles in the
flow. Moreover,  
one has  that $\FF=(\nabla\bm\xi)^{-1}$.

As $\FF$ depends on $\xx$,  equalities
\eqref{ultimate}--\eqref{transport-xi} will be assumed to  
hold for  almost all  $\xx  \in \yy(t,\Omega)$.   
Here,  we take advantage of   
the boundary condition $\vv{\cdot}\nn=0$ imposed below,  with $\nn$ being
the outward unit normal to the boundary of the actual domain. This boundary
condition  implies that   
the actual domain $\varOmega  = \yy(t,\Omega)$  does not evolve in
time, i.e., $\varOmega=\Omega$.
 The stress $\TT$ and the pressure $p$ defined in
\eqref{Euler1-referential-selfgravit} and \eqref{Euler1-fluid-selfgravit}
are also functions of $\xx$. In particular, all  models  under
consideration in the following will be
fully Eulerian.

\section{Self-gravitating inhomogeneous viscoelastic media}\label{sec-single}

We  consider  a simple model  for a   self-gravitating bounded body with
a fixed shape.  
Planets and moons are typically composed  by  many
components.  By referring to the Earth, as well as to other
planets and moons of the solar system, one should minimally consider
three components, namely, metals, silicates, and  a  gaseous
atmosphere.  
 The latter is often adopted for  the so-called {\it sticky-air approach}
(cf.\ \cite{CSGD12CNST}) to  allow for a fixed  and smooth
 domain $\varOmega$.  
We are thus led to consider
here a spatially inhomogeneous material with a given  referential  mass density
$ \rho_{\rm ref} = \rho_{\rm ref}   (\XX)>0$ at  some initial  time $t=0$
and stored energy $\varphi_{\rm ref}(\XX,\cdot):\R^{3\times3}\to\R$ acting
on $\FF$. 
We  preliminarily assume that the body follows  the {\it Kelvin-Voigt
  rheology}, which is the simplest viscoelastic solid-type rheology.   In
the following, two inhomogeneous viscosity coefficients $\nu_1$ and $\nu_2$
will be introduced.

We will use the shorthand notation $(\cdot)^{\bm\xi}$  in order to
indicate the  composition
with $\bm\xi$  within a referential quantify when substituting  
$\XX$  by  $\bm\xi(t,\xx)$.  In particular, 
$\FF=\FF_\RRR^{\bm\xi}$  and  $\vv=\vv_\RRR^{\bm\xi}$. This
notation is  meant  to
emphasize the spatial inhomogeneity of the material.

The referential density $\rhoR=\rhoR(\XX)$ in the deformed configuration
is $\rhoR=\rho_{\rm ref}/\!\det\FF_\RRR$ while the actual Eulerian mass
density $\varrho=\varrho(t,\xx)$ is
$\rhoRxi=\rho_{\rm ref}^{\bm\xi}/\!\det\FF_\RRR^{\bm\xi}$, i.e.,
\begin{align}\label{rho=rho0/detF}
\varrho=\frac{\rho_{\rm ref}^{\bm\xi}}{\det\FF}=\det(\nabla\bm\xi)\rho_{\rm ref}^{\bm\xi}\,.
\end{align} 
Likewise, the referential stored energy  reads 
$\varphi_\RRR^{}=\varphi_\RRR^{}(\XX,\FF_\RRR)=\varphi_{\rm ref}^{}(\XX,\FF_\RRR)/\!\det\FF_\RRR$
and the actual Eulerian stored energy  is  $\varphi=\varphi(t,\xx,\FF)
=\varphi_{\rm ref}^{\bm\xi}(t,\xx,\FF)/\!\det\FF$. Note that now
$\varphi$  is dependent on time $t$.   
In what follows, we will often  simply   write
$\varphi_{\rm ref}^{\bm\xi}(\FF)$ instead of
$\varphi_{\rm ref}^{\bm\xi}(t,\xx,\FF)$.

 Given equation  \eqref{ultimate},  
relation \eqref{rho=rho0/detF} is
equivalent to the evolution-and-transport equation for the actual mass density
\begin{align}\label{cont-eq+}
\DT\varrho=-\varrho\,{\rm div}\,\vv\ \ \ \text{ with the initial condition}\ \
\varrho|_{t=0}^{}=\frac{\rho_{\rm ref}^{\bm\xi}|_{t=0}}{\det\FF|_{t=0}^{}}\,.
\end{align}

 A distinct advantage of working in a fully Eulerian setting is
the simplicity with which one can incorporate  in the model
interactions with actual fields. Relevant to our endeavor is in particular the
gravitational acceleration field  $\gg=-\nabla\GRAVPOT$ 
ensuing from the {\it gravitational potential} $\GRAVPOT$.
\DELETE{We  introduce a simplification of the model by posing the system
in a bounded and smooth albeit possibly large domain $U\subset \R^3$,
which plays the role of {\it universe}. In addition, we assume that the
gravitational potential $V$ takes the constant value $V_\text{\sc b}^{}$
outside $U$. The rationale of this choice is
that, by possibly taking $U$ large enough, the effect  of assuming
the universe  $U$ to be bounded  on the dynamic
of the planet   is expected to be negligible. Note that our system will be
completely independent of the constant $V_\text{\sc b}^{}$. This is
particularly relevant in connection with the choice of $U$. Note
that it is well-known that the actual value of the gravitational
potential on the Earth surface depends on whether one considers only
a single planet, the whole Solar System, or the whole galaxy Milky
Way,  and is given  roughly by 60 MJ/kg, 900 MJ/kg, or more than 130 GJ/kg,
respectively. Fixing a specific value  $V_\text{\sc b}^{}$ in \eqref{BC}
will in fact be shown to be immaterial in our framework,  see
\eqref{calculus-selfgravit} below.} 

The gravitational potential $\GRAVPOT$ is governed by the Poisson equation
\begin{align}\label{Poisson}
\Delta\GRAVPOT=\begin{cases}\TTT4\uppi\EEE\GRAVCONST\varrho
\ \ \text{ with }\ \varrho=\det(\nabla\bm\xi)\rho_{\rm ref}^{\bm\xi}
&\text{on }\ \varOmega\,,\\
\TTT4\uppi\EEE\GRAVCONST\varrho_{\rm ext}^{}
&\text{on }\ \Universe\setminus\varOmega,
\end{cases}
\end{align}
where $\GRAVCONST$ is the gravitational constant
\TTT ($\doteq 6.674{\times}10^{-11}$m$^3$kg$^{-1}$s$^{-2}$) \EEE and
$\varrho_{\rm ext}^{}=\varrho_{\rm ext}^{}(t,\xx)$  is  a given
external,  possibly time-dependent  mass density distributed around
$\varOmega$,  which may model tidal effects \TTT caused by moons orbiting
the planet\EEE.  In fact,
we will introduce a further simplification by assuming that the mass
$\varrho_{\rm ext}^{}$ does not feel the presence of the mass
$\varrho$. 
In what follows, we will consider both %
$\varrho$ and $\varrho_{\rm ext}^{}$  to be defined on the whole
$\Universe$ by extending them to $0$ outside   
$\varOmega$
and $\Universe\setminus\varOmega$,  respectively.   

The  geometric setting of the model   
is illustrated in Figure~\ref{fig1}.
\begin{center}
\begin{my-picture}{6}{6}{fig1}
\psfrag{W}{$\varOmega$}
\psfrag{U}{ $\Universe$}
\psfrag{V=0}{\small $\GRAVPOT=\GRAVPOT_\text{\sc b}^{}$}
\psfrag{v=0}{\small $\vv{\cdot}\nn=\bm0$} 
\psfrag{r,v,xi,J}{\small $\varrho,\,\vv,\,\bm\xi,\,J$}
\psfrag{rho-ext}{\small $\varrho_{\rm ext}^{}$}
\psfrag{flux=0}{\small $\nn{\cdot}\nabla\bm\mu=\bm0$} 
\hspace*{-2em}\includegraphics[width=20em]{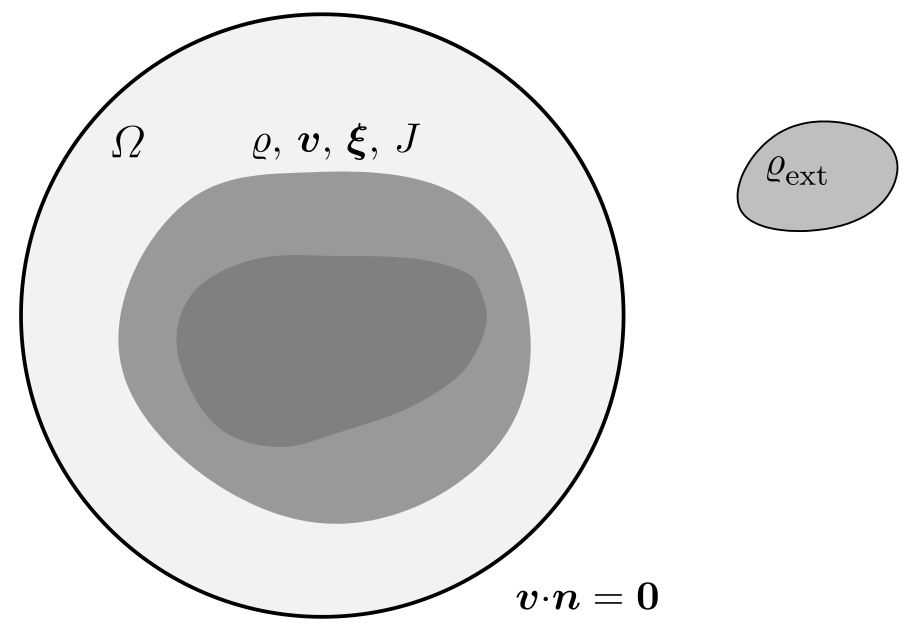}
\end{my-picture}
\nopagebreak
\\
{\small\sl\hspace*{-3em}Fig.\,\ref{fig1}:~\begin{minipage}[t]{34em}
The schematic geometry of the fixed bounded domain $\varOmega\subset\Universe$
where the self-gravitating medium
 is evolving with velocity $\vv$ \TTT whereas \EEE
the external mass density $\varrho_{\rm ext}^{}=\varrho_{\rm ext}^{}(t,\xx)$ outside
$\varOmega$ is prescribed.
The boundary conditions on $\pl\varOmega$
are depicted  too,  with 
$\nn$ denoting the 
unit normal to the boundary $\pl\varOmega$.
\end{minipage}
}
\end{center}

\medskip

The conservative and the dissipative parts of the Cauchy stress 
acting on the body are defined  
standardly (cf.~Remark~\ref{rem-grad}) as
\begin{subequations}
\label{T+D}
\begin{align}
&\label{T+D1}\TT=\frac{\big[\varphi_{\rm ref}^{\bm\xi}\big]_\FF'(\FF)\FF^\top\!\!\!}{\det\FF}=
\big[\varphi_{\rm ref}^{\bm\xi}\big]_\FF'((\nabla\bm\xi)^{-1}){\rm Cof}(\nabla\bm\xi)
\ \ \text{ and }\ \
\\\label{T+D2}&
\DD=\nu_1^{\bm\xi}\EE(\vv)-{\rm div}\big(\nu_2^{\bm\xi} 
   |\nabla\EE(\vv)|^{q-2}\nabla\EE(\vv)\big)\,.
\end{align}
\end{subequations} 
 The Kelvin-Voigt rheology  corresponds to  
considering the total stress as $\TT+\DD$.  Note that $\DD$ contains a
standard Newtonian term $\nu_1^{\bm\xi}\EE(\vv)$ and the hyperviscous stress 
$-{\rm div}\big(\nu_2^{\bm\xi} |\nabla\EE(\vv)|^{q-2}\nabla\EE(\vv)\big)$, see
Remark \ref{rem-grad} for a discussion.

The full model system is formulated in the unknowns $(\vv,\bm\xi,\GRAVPOT)$
and results from the momentum equilibrium, \eqref{ultimate}, namely,  
\begin{subequations}
  \label{Euler-referential-selfgravit}
 \begin{align}
&&&\varrho\DT\vv={\rm div}\big(\TT{+}\DD\big)
      -  \varrho \nabla\GRAVPOT\!\!\!\!
&&\text{ on }\ \varOmega\ \ \text{ with $\TT$ and $\DD$ from \eqref{T+D}
and $\varrho$ from \eqref{rho=rho0/detF}},&&&&
\label{Euler1-referential-selfgravit}
\\\label{Euler3-referential-selfgravit}
&&&
\DT{\bm\xi}=\bm0&&\text{ on }\ \varOmega\,,
\\\label{Euler2-referential-selfgravit}
&&&\Delta\GRAVPOT=\TTT4\uppi\EEE\GRAVCONST\big(
\varrho{+}\varrho_{\rm ext}^{}\big)&&\text{ on }\ \Universe\,.
\end{align}\end{subequations}
We  complement this system by  the boundary conditions
\begin{align}
&\label{BC}\vv{\cdot}\nn=0\,,\ \ \ 
\big((\TT{+}\DD)\nn-\divS(
 \nu_2^{\bm \xi}|\nabla\EE(\vv)|^{q-2}\nabla\EE(\vv)
\nn)\big)_\text{\sc t}^{}=\bm0\,,\ \ \ \nabla\EE(\vv){:}(\nn{\otimes}\nn)=\bm0\ \
\text{ on }\pl\varOmega.
\end{align}
Here, $(\cdot)_\text{\sc t}^{}$ indicates  the tangential
component of a vector  at $\partial\varOmega$\DELETE{ and
$\GRAVPOT_\text{\sc b}^{}$ is  the above mentioned arbitrary
constant for the value of the gravitational field outside $U$}. 
The  $(d{-}1)$-dimensional surface divergence is
defined as
\begin{align}\label{def-divS}
\divS={\rm tr}(\nablaS)\ \ \ \text{ with }\ \
\nablaS \bulet=\nabla \bulet-\frac{\partial \bulet}{\partial\nn}\nn\,,
\end{align}
where ${\rm tr}(\cdot)$ is the trace of a
$(d{-}1){\times}(d{-}1)$-matrix and
$\nablaS  $ denotes the surface gradient. 
Let us again remark the crucial role of the impenetrability boundary
condition $\vv{\cdot}\nn=0$, indeed allowing system
(\ref{Euler-referential-selfgravit}a,b) to be formulated in the fixed domain
$\varOmega$.
\TTT
The gravitational equation \eq{Euler2-referential-selfgravit}
implicitly requires to specify a value for the potential at
$\infty$. The actual choice of this value is immaterial, as the
coupling is realized via the 
gravitational accelleration only.  A customary choice is to set $\lim_{|\xx|\to \infty}\GRAVPOT(\xx)=0$, so that the
following representation formula holds
\begin{align}\label{V-from-integral}
\GRAVPOT(t,\xx)=-\GRAVCONST\!\int_{\R^3}\!\!
\frac{\varrho(t,\wt\xx)+\varrho_{\rm ext}^{}(t,\wt\xx)}
{|\xx{-}\wt\xx|}\,\d\wt\xx\,.
\end{align}
The gravitational acceleration $-\nabla\GRAVPOT$ in
\eq{Euler1-referential-selfgravit} is
\begin{align}\label{nabla-V}
-\nabla\GRAVPOT(t,\xx)=\GRAVCONST\!\int_{\R^3}\!\!
\frac{(\varrho(t,\wt\xx){+}\varrho_{\rm ext}^{}(t,\wt\xx))
(\wt\xx{-}\xx)}{|\wt\xx{-}\xx|^3}\,\d\wt\xx\,
\end{align}
and the gravitational force $-\varrho(\xx)\nabla\GRAVPOT(\xx)$ in
\eq{Euler1-referential-selfgravit} is then
$$
-\varrho(t,\xx)\nabla\GRAVPOT(t,\xx)=
\int_\varOmega\frac{\varrho(t,\xx)
\varrho(t,\wt\xx)(\wt\xx{-}\xx)}{|\wt\xx{-}\xx|^3}\,\d\wt\xx
+\varrho(t,\xx)\int_{\R^3\setminus\varOmega}\!\!\!\!\!
\frac{\varrho_{\rm ext}^{}(t,\wt\xx)
(\wt\xx{-}\xx)}{|\wt\xx{-}\xx|^3}\,\d\wt\xx
\ \ \text{ for}\ \xx\in\varOmega\,.
$$
The first right-hand-sided term corresponds to Newton's universal law of 
gravitation. By substitution into \eq{Euler1-referential-selfgravit},
this, in fact, allows for elimination of $\GRAVPOT$ from the reduced system
(\ref{Euler-referential-selfgravit}a,b).

\AAA 
Subsequently, we will use that the dependence of $V=\mathsf V[\wh\varrho]$ with
$\wh\varrho=\varrho+ \varrho_{\rm ext}$ is linear, and for all open bounded sets $U
\subset \R^3$ the following operators are bounded: 
\begin{equation}
  \label{eq:V.lin.oper}
  \mathsf V :\left\{ \begin{array}{ccc} L^2(U) &\to & H^2(U); \\ 
    \wh\varrho &\mapsto &V=\mathsf V[\wh\varrho]. 
  \end{array} \right. 
 \quad \text{and} \quad
  \nabla \mathsf V :\left\{ \begin{array}{ccc} L^2(U) &\to & H^1(\R^3;\R^3); \\ 
    \wh\varrho &\mapsto &\nabla V=\nabla\mathsf V[\wh\varrho]. 
  \end{array} \right.
\end{equation}
Note that $V$, $\nabla V$, $\nabla^2 V$ decay like $1/|\xx|$, $1/|\xx|^2$, and
$1/|\xx|^3$, respectively, hence $V\not\in L^2(\R^3)$ but $\nabla V\in
L^2(\R^3;\R^3) $   and $ \nabla^2 V\in L^2(\R^3;\R^{3 \times 3})$. 
\EEE

The energetics  of system \eqref{Euler-referential-selfgravit} can
be revealed  by testing \eqref{Euler1-referential-selfgravit} by $\vv$
and \eqref{Euler2-referential-selfgravit} by $\pdt{}\GRAVPOT$.
The former test quite  standardly employs the Green formula over $\varOmega$
(for the $\nu_2$ term  one uses the formula twice, combined with a surface
Green formula on $\pl\varOmega$) in view of the first and the second boundary
conditions in \eqref{BC}. In fact, following
\cite[Sec.\,3]{Roub22IFSV} we start from \eqref{T+D1} and compute
\begin{align}
  \TT&=\frac{\big[\varphi_{\rm ref}^{\bm\xi}\big]_\FF'(\FF)\FF^\top\!\!\!}{\det\FF}
  = \frac{\big[\varphi_{\rm ref}^{\bm\xi}\big]_\FF'(\FF)-
  \varphi_{\rm ref}^{\bm\xi}(\FF)\FF^{-\top}\!\!\!}{\det\FF}\FF^\top\!
  + \frac{\varphi_{\rm ref}^{\bm\xi}(\FF)}{\det\FF}\bbI\nonumber\\
&=\bigg( \frac{\big[\varphi_{\rm ref}^{\bm\xi}\big]_\FF'(\FF)}{\det\FF}
- \frac{\varphi_{\rm ref}^{\bm\xi}(\FF){\rm Cof}\FF}{(\det\FF)^2}\bigg)\FF^\top\!
+ \frac{\varphi_{\rm ref}^{\bm\xi}(\FF)}{\det\FF}\bbI
= \bigg[ \frac{\varphi_{\rm ref}^{\bm\xi}(\FF)}{\det\FF}\bigg]'_\FF\FF^\top\!\!
+ \frac{\varphi_{\rm ref}^{\bm\xi}(\FF)}{\det\FF}\bbI\,,\label{Told}
\end{align}
where $\bbI$ denotes the identity matrix. Moreover, we have that 
\begin{subequations}
  \begin{align}
\frac{\partial}{\partial t}\bigg(\frac{\varphi_{\rm ref}^{\bm\xi}(\FF)}{\det\FF}\bigg)
  & = \frac{\big[[\varphi_{\rm ref}]'_\XX\big]^{\bm\xi}\!}{\det\FF}{\,\cdot}
  \frac{\partial\bm\xi}{\partial t} +
  \bigg[\frac{\varphi_{\rm ref}(\FF)}{\det\FF}\bigg]'_\FF{:} \frac{\partial
                                          \FF}{\partial t}\ \ \text{ and}\\
\nabla\bigg(\frac{\varphi_{\rm ref}^{\bm\xi}(\FF)}{\det\FF}\bigg)
\SSS{\cdot}\vv\EEE
 & = \frac{\big[[\varphi_{\rm ref}]'_\XX\big]^{\bm\xi}\!}{\det\FF}\,{\cdot}
 (\vv{\cdot} \nabla){\bm \xi}+ \bigg[\frac{\varphi_{\rm ref}(\FF)}{\det\FF}\bigg]'_\FF{:} (\vv{\cdot} \nabla)\FF.
  \end{align}\label{for1}
\end{subequations}

By testing the force $-{\rm div}(\TT{+}\DD)$ by the velocity $\vv$ and
using the boundary conditions  \eqref{BC} we hence get
\begin{align*}
 &\hspace{-4em} -\int_\varOmega\!{\rm div}(\TT{+}\DD){\cdot} \vv \, \d
  \xx\stackrel{\eqref{BC}}{=} \int_{\varOmega} (\TT{+}\DD){:}\nabla \vv
  \, \d \xx
 \\&\quad
\stackrel{\eqref{Told}}{=}\int_\varOmega\bigg[\frac{\varphi_{\rm
  ref}(\FF)}{\det\FF}  \bigg]'_\FF\FF^{\top} {:} \nabla \vv + \frac{\varphi_{\rm
  ref} (\FF)\bbI}{\det\FF}{:} \nabla \vv + \DD{:}\nabla
  \vv \, \d \xx
\\
&\quad\ =\int_\varOmega\bigg[\frac{\varphi_{\rm ref}(\FF)}{\det\FF}\bigg]'_\FF{:}(\nabla\vv)\FF + \frac{\varphi_{\rm
  ref} (\FF)}{\det\FF}{\rm div}\,\vv + \DD{:}\nabla
  \vv \, \d \xx
\\
&\quad\stackrel{\eqref{ultimate}}{=}\int_\varOmega\bigg[\frac{\varphi_{\rm
  ref}(\FF)}{\det\FF}  \bigg]'_\FF{:} \DT \FF + \frac{\varphi_{\rm
  ref} (\FF)}{\det\FF}{\rm div}\,\vv+ \DD{:}\nabla
  \vv \, \d \xx\\
&\quad \stackrel{\eqref{for1}}{=}\int_{\varOmega}\frac{\partial}{\partial t} \bigg( \frac{\varphi_{\rm
  ref}^{\bm\xi}(\FF)}{\det\FF}\bigg) +
  \nabla \bigg( \frac{\varphi_{\rm
  ref}^{\bm\xi}(\FF)}{\det\FF}\bigg){\cdot} \vv +\frac{\varphi_{\rm
  ref} (\FF)}{\det\FF}{\rm div}\,\vv\, \d \xx
\\&\hspace{10em}
- \int_{\varOmega}
  \frac{\big[[\varphi_{\rm
  ref}]'_\XX\big]^{\bm\xi}\!}{\det\FF}\,{:}
  \bigg(\frac{\partial\bm\xi}{\partial t}+ \vv {\cdot} \nabla {\bm\xi} \bigg) + \DD{:}\nabla
  \vv \, \d \xx\\
&\quad \stackrel{\eqref{transport-xi}}{=} \frac{\d}{\d
  t}\int_{\varOmega} \frac{\varphi_{\rm
  ref}^{\bm\xi}(\FF)}{\det\FF}\, \d \xx + \int_{\varOmega}
  {\rm div} \bigg( \frac{\varphi_{\rm
  ref}^{\bm\xi}(\FF)}{\det\FF}\, \vv\bigg) \, \d \xx + \int_\varOmega\DD{:}\nabla
  \vv \, \d \xx\\
&\quad\ 
=\frac{\d}{\d t}\int_{\varOmega} \frac{\varphi_{\rm
  ref}^{\bm\xi}(\FF)}{\det\FF}\, \d \xx + \int_{\partial \varOmega}
   \frac{\varphi_{\rm
  ref}^{\bm\xi}(\FF)}{\det\FF}\, ( \vv{\cdot}\nn) \, \d S + \int_\varOmega\DD{:}\nabla
  \vv \, \d \xx\\
&\quad\ = \frac{\d}{\d
  t}\int_{\varOmega} \frac{\varphi_{\rm
  ref}^{\bm\xi}(\FF)}{\det\FF}\, \d \xx + \int_\varOmega
  \nu_1^{\bm\xi}|\EE(\vv)|^2 +
  \nu_2^{\bm\xi}|\nabla \EE(\vv)|^q \d \xx.
\end{align*}

By testing the inertial force $\varrho\DT\vv$ by $\vv$, it is paramount to
use \eqref{Euler3-referential-selfgravit} together with 
\eqref{rho=rho0/detF}  and  $\FF=(\nabla\bm\xi)^{-1}$, 
inducing  
the
continuity equation \eqref{cont-eq+}.  In particular, one has
$\pdt{}(\varrho|\vv|^2/2)=\varrho\vv{\cdot}\pdt{}\vv-\frac12{\rm div}(\varrho\vv)|\vv|^2$.
Hence, by using again the Green formula and $\vv{\cdot}\nn=0$, we obtain
\begin{equation}
\frac{\d}{\d t}\int_\varOmega\frac{\varrho|\vv|^2}{2} \,\d\xx
=\int_\varOmega\varrho\DT\vv{\cdot}\vv\,\d\xx.\label{eq:referlater}
\end{equation}
On the other hand, testing \eqref{Euler2-referential-selfgravit} on
$\pdt{}\GRAVPOT$, we obtain  
 \begin{align}\nonumber
&\frac{\d}{\d t}\int_{\Universe}\frac{|\nabla\GRAVPOT|^2}{\TTT8\uppi\EEE\GRAVCONST}\,\d\xx
=
\int_{\Universe}
\frac{\Delta\GRAVPOT}{\TTT4\uppi\EEE\GRAVCONST}\pdt{\GRAVPOT}\,\d\xx
   \stackrel{\eqref{Euler2-referential-selfgravit}}{=}
-\!\int_{\Universe}(\varrho{+}\varrho_{\rm ext}^{})\pdt{\GRAVPOT}\,\d\xx\nonumber
\nonumber\\&\qquad\ \nonumber 
=\int_{\Universe}\bigg(\pdt\varrho{+}\pdt{\varrho_{\rm ext}}\bigg)\GRAVPOT\,\d\xx
-\frac{\d}{\d t}\int_{\Universe}(\varrho{+}\varrho_{\rm ext}^{})\GRAVPOT\,\d\xx
\\&\qquad\nonumber 
   \stackrel{\eqref{cont-eq+}}{=}
   -\int_\varOmega {\rm div}(\varrho\vv)\, \GRAVPOT\,\d\xx
   +\int_{\Universe{\setminus}\varOmega}\!\!\pdt{\varrho_{\rm ext}}\GRAVPOT\,\d\xx
-\frac{\d}{\d t}\int_{\Universe}(\varrho{+}\varrho_{\rm ext}^{})\GRAVPOT\,\d\xx
   \\&\qquad\ 
=\int_\varOmega\varrho\vv{\cdot}\nabla \GRAVPOT \,\d\xx
+\int_{\Universe{\setminus}\varOmega}\!\!\pdt{\varrho_{\rm ext}}\GRAVPOT \,\d\xx
-\frac{\d}{\d t}\int_\varOmega\varrho \GRAVPOT \,\d\xx
-\frac{\d}{\d t}\int_{\Universe{\setminus}\varOmega}\varrho_{\rm ext}^{} \GRAVPOT
\,\d\xx\, 
\label{calculus-selfgravit}\end{align}
with $\varrho$ from \eqref{rho=rho0/detF},  where we also used the
fact that $\int_{\varOmega}\varrho\, \d\xx$ is constant. 
\DELETE{Note that  the Green formula  on  $\Universe$ used in
\eqref{calculus-selfgravit} hinges  on the boundary condition $\GRAVPOT=
\GRAVPOT_\text{\sc   b}^{} $ on $\pl\Universe$  while the Green
formula on $\varOmega$ used again the boundary condition
$\vv{\cdot}\nn=0$.  The equality in \eqref{calculus-selfgravit} is
nonetheless independent of the value $\GRAVPOT_\text{\sc b}^{} $.}
 By noticing  that the term
$\varrho\vv{\cdot}\nabla\GRAVPOT=\varrho\nabla\GRAVPOT{\cdot}\vv$
arises also when testing \eqref{Euler1-referential-selfgravit} by $\vv$ where
it is to be substituted from \eqref{calculus-selfgravit}, we obtain the
mechanical energy-dissipation balance
\begin{align}\nonumber
&\hspace*{0em}\frac{\d}{\d t}\Bigg(\int_\varOmega\!\!\!
\linesunder{\frac{\varrho}2|\vv|^2}{kinetic}{energy}\!\!\!\!+\!\!\!
  \linesunder{\frac{\varphi_{\rm ref}^{\bm\xi}(\FF)}{\det\FF}}{stored}{energy}
\!\!+\!\!\!\!\!\!\!\!\!\!\morelinesunder{
\varrho\GRAVPOT
}{\ \ energy of $\varrho$ in}{\ \ gravitational}{field $\GRAVPOT$}\!\!\!\!\!\!\!\!\!\!\d\xx
\ \ +\ \int_{\Universe}\!\!\!\!\!\!\!\!\morelinesunder{\frac{|\nabla\GRAVPOT|^2\!}{\TTT8\uppi\EEE\GRAVCONST}}{energy of}{gravitational}{field}\!\!\!\!\!\!\d\xx
+\!\!\int_{\Universe\setminus\varOmega}
\!\!\!\!\!\!\!\!\morelinesunder{\varrho_{\rm ext}^{}\GRAVPOT}{energy of $\varrho_{\rm ext}^{}$}{\ \ \ in gravitational}{field $\GRAVPOT$}\!\!\!\!\!\!\!\!\!\!\!\!\d\xx\Bigg)
\\&\qquad\qquad\qquad\qquad\qquad
+\!\int_\varOmega\!\!\linesunder{\nu_1^{\bm\xi}|\EE(\vv)|^2
+\nu_2^{\bm\xi}|\Nabla\EE(\vv)|^{ q }_{_{_{}}}\!}{$=:\xi$ dissipation rate due}{to viscosity}\d\xx
=\int_{\Universe\setminus\varOmega}\!\!\!\!\!\!\!\!\!\!\!\!\morelinesunder{\pdt{\varrho_{\rm ext}}\GRAVPOT}{power of}{external mass in}{gravitational field}\!\!\!\!\!\!\!\!\!\d\xx\,.
\label{mech-engr-selfgravit}
\end{align}

\begin{remark}[{\sl Variational structure of
\eqref{Euler-referential-selfgravit}}]\upshape
To  elucidate the  variational structure  of  the system
\eqref{Euler-referential-selfgravit}  requires some care, because the
potential equation \eqref{Poisson} for $\GRAVPOT$ provides a concave
contribution to the free energy
\begin{align*}
(\yy,\GRAVPOT)\mapsto
\int_\Omega\varphi_{\rm ref}(\XX,\nabla\yy)
+\varrho_{\rm ref}\big(\XX,\GRAVPOT{\circ}\yy\big)\,\d\XX
+\int_{\Universe\setminus\varOmega}\!\!\!\varrho_{\rm ext}^{}(\xx)\GRAVPOT(\xx)\,\d\xx
-\int_{\Universe}\frac{|\nabla\GRAVPOT(\xx)|^2}{\TTT8\uppi\EEE\GRAVCONST}\,\d\xx\,.
\end{align*}
Note that the latter features a mixture of referential and
actual terms. By taking the variation of the free energy with
respect to the gravitational potential $\GRAVPOT$ and using
$\GRAVPOT{\circ}\yy(\XX)=\GRAVPOT(\xx)$ we obtain \eqref{Poisson}. On
the other hand, the variation with respect to $\yy$ produces the first
Piola-Kirchhoff stress tensor and the gravitational
force $\rhoR(\nabla\GRAVPOT{\circ}\yy)$, which correspond to $\TT$ and
$\rhoR\nabla\GRAVPOT=\varrho\nabla\GRAVPOT$ when written in the
Eulerian  setting of  \eqref{Euler1-referential-selfgravit}. 
\end{remark}

\begin{remark}[{\sl Gradient theories in rates}]
\label{rem-grad}\upshape
 Higher-order theories in solid mechanics are well established and
used for various reasons.  By introducing a further length scale into the
problem, additional hyperstresses occur, which in turn usually contribute
crucial compactness for the mathematical analysis. Such high-order models are
generally referred to as {\it nonsimple materials}. Both the conservative and
the dissipative stress can feature higher gradients. In
\eqref{Euler-referential-selfgravit} we consider an hyperstress on the dissipative
part, an option which is particularly well tailored to rate formulations,
having the advantage  to provide  additional regularity for the velocity
field $\vv$. Our approach 
follows the theory by E.~Fried and M.~Gurtin \cite{FriGur06TBBC}, as already
 considered  in the general nonlinear context of {\it multipolar
  fluids} by J.~Ne\v cas at
al.\ 
\cite{Neca94TMF,NeNoSi91GSCI,NecRuz92GSIV,NecSil91MVF}
 and as originally  inspired by
R.~A.~Toupin \cite{Toup62EMCS} and R.~D.~Mindlin \cite{Mind64MSLE}.
\end{remark}

\begin{remark}[{\sl  Anisothermal  extension}]\label{rem-anisotherm0}\upshape
  Heat  exchange and 
  transfer  plays an important role with respect to the differentiation
  phenomenon in self-gravitating planets and moons. Although presently
neglected in our model, thermal effects could also be considered. On the one
  hand, one may let material parameters be dependent on temperature, here
  denoted by $\theta$.  On the other hand, by
  \TTT making \EEE the  
  (referential) free energy
$\psi_{\rm ref}=\psi_{\rm ref}(\XX,\FF,\theta)$ \DELETE{to be}
additively decomposed as 
$\psi_{\rm ref}(\XX,\FF,\theta)=\varphi_{\rm ref}(\XX,\FF)
+\COUPLING_{\rm ref}(\XX,\FF,\theta)$ with $\gamma_{\rm ref}(\XX,\FF,0)=0$, the
heat equation reads 
\begin{align}\label{heat-eq}
\pdt\W+{\rm div}\big(\vv\W{+}\jj\big)=\xi
+\frac{[\COUPLING_{\rm ref}^{\rm\xi}]'_{\FF}(\FF,\theta)}{\det\FF}
{:}\DT\FF
\ \ \ \text{ with }\ \ \W=\frac{\COUPLING_{\rm ref}^{\rm\xi}(\FF,\theta)
{-}\theta[\COUPLING_{\rm ref}^{\rm\xi}]_\theta'(\FF,\theta)}
{\det\FF}\,,
\end{align}
where the heat production rate $\xi$  is specified in
\eqref{mech-engr-selfgravit} and  
the heat flux $\jj$ is governed by the Fourier law
$\jj=-k^{ \boldsymbol\xi}\nabla\theta$ with a heat-conduction coefficient
$k=k(\XX,\det\FF,\theta)$.   Then,  the viscosity coefficients 
$\nu_1$ and $\nu_2$  can be assumed to depend  on temperature,  as
well.  The physical meaning of  the quantity $\W$ is of that of 
the thermal part of the (actual) internal energy.  By complementing
relation \eqref{heat-eq} to the system \eqref{Euler-referential-selfgravit}, the
ensuing  anisothermal  coupled model  then reproduces the expected
energetics. In addition, it complies with the  Clausius-Duhem entropy
inequality, cf.\ \cite{Roub22TVSE}  and it is hence 
thermodynamically consistent.
\end{remark}

\section{Kelvin-Voigt/{}Navier-{}Stokes viscoelastic fluid}\label{sec-fluids}
 Differentiation by self-gravitation in planets and moons occurs
on  very long time scales.  Within  such a  time scale, solid-type
rheologies (as the Kelvin-Voigt one from Section~\ref{sec-single})
have limited relevance and one should resort to fluid rheologies in
the deviatoric components instead. Different options are
available. Classical choices in geophysics are the Andrade and the
Jeffreys rheologies, both allowing  for  the propagation of elastic
shear waves. Such waves are however of  relatively low  importance
on the time scale of planetary evolution. As such, we resort here to a
simpler variant, which goes under the name of Newton or Stokes or, in
the current convective setting, of Navier-Stokes rheology and 
 does not allow for propagation of shear waves. 
 The Navier-Stokes rheology can be obtained in the frame of the
above introduced  
Kelvin-Voigt model  by assuming that  
the  elastic shear  response  vanishes,
i.e.,  that  $\varphi_{\rm ref}$ depends only on the 
volumetric 
part $(\det\FF)^{1/3}\bbI$
of $\FF$.
Thus,  we let 
\begin{align}
\varphi_{\rm ref}(\XX,\FF)=\phi_{\rm ref}(\XX,J)\ \ \ \text{  with  }\ \ \ J=\det\FF
\label{stored-fluid}
\end{align}
for some $\phi_{\rm ref}:\varOmega\times(0,\infty)\to\R$, cf.,\ e.g., \cite[p.\,10]{MarHug83MFE} or also \cite{Roub22IFSV}.
By  recalling that  $\det'(\cdot)={\rm Cof}(\cdot)$ and
$F^{-1}={\rm Cof}F^\top\!/\!\det F$, the conservative part
of the Cauchy stress reduces to 
\begin{align}\nonumber
\TT&=\frac{[\varphi_{\rm ref}]_\FF'(\XX,\FF)\FF^\top\!\!}{\det\FF}
 =[\phi_{\rm ref}]_J'(\XX,\det\FF)\frac{\det'(\FF)\FF^\top\!\!}{\det\FF}
\\&=[\phi_{\rm ref}]_J'(\XX,\det\FF)\frac{({\rm Cof}\,\FF)\FF^\top\!\!}{\det\FF}
  =[\phi_{\rm ref}]_J'(\XX,J)\bbI\,,
\label{T-p}\end{align}
 (compare with \eqref{Told}), 
where $[\phi_{\rm ref}]_J'$ has a physical interpretation  of  a (negative) pressure.

 Under  the fluidic ansatz \eqref{stored-fluid},  by taking
\eqref{T-p} and \eqref{rho=rho0/detF} into account, 
system
\eqref{Euler-referential-selfgravit} reads
\begin{subequations}\label{Euler-fluid-selfgravit}\begin{align}
\nonumber
&\varrho\DT\vv=
{
\rm div}\big(
\nu_1^{\bm\xi}\EE(\vv)-{\rm div}\big(\nu_2^{\bm\xi}|\nabla\EE(\vv)|^{q-2}\nabla\EE(\vv)\big)\big)
-\nabla p
       - \varrho\nabla\GRAVPOT\ \ \ 
 \\& \hspace{8.5em}\text{ with }\ \
p=-[\phi_{\rm ref}^{\bm\xi}]_J'\Big(\frac1{\det(\nabla\bm\xi)}\Big)
\ \ \text{ and }\ \ \varrho=\det(\nabla\bm\xi)\rho_{\rm ref}^{\bm\xi}
\,,
\label{Euler1-fluid-selfgravit}
\\[-.5em]\label{Euler3-fluid-selfgravit}
&
\DT{\bm\xi}=\bm0,\hspace{6em}\text{on }\ \varOmega\,, 
\\\label{Euler2-fluid-selfgravit}
&\Delta\GRAVPOT=\TTT4\uppi\EEE\GRAVCONST\big(
\varrho{+}\varrho_{\rm ext}^{}\big)\ \ \
\text{ on }\ \Universe\,.
\end{align}\end{subequations}
 System 
\eqref{Euler-fluid-selfgravit} is  usually  referred to as 
 the  
compressible Navier-Stokes-Poisson  system. 
We  complement it   with the boundary conditions \eqref{BC}.

The energetics  for system \eqref{Euler-fluid-selfgravit} follows along the
same lines of that of system  
\eqref{Euler-referential-selfgravit}. Let us explicitly show  how to handle
the conservative part of the stress, i.e.,  the term 
${\rm div}\,\TT{\cdot}\vv$, which now reads as $\nabla p{\cdot}\vv$ with
$p=-[\phi_{\rm ref}^{\bm\xi}]_J'(J)$.  The analogous of computations
\eqref{for1} in the fluidic setting \eqref{stored-fluid} are 
\begin{subequations}\label{calculus}\begin{align}
&\pdt{}\Big(\frac{\phi_{\rm ref}^{\bm\xi}(J)}{J}\Big)=
\frac{\big[[\phi_{\rm ref}]_\XX'\big]^{\bm\xi}(J)}{J}{\cdot}\pdt{\bm\xi}+
\Big[\frac{\phi_{\rm ref}^{\bm\xi}(J)}J\Big]_J'
\pdt J
\\&\nabla\Big(\frac{\phi_{\rm ref}^{\bm\xi}(J)}{J}\Big){\cdot}\vv=
\frac{\big[[\phi_{\rm ref}]_\XX'\big]^{\bm\xi}(J)}{J}{\cdot}(\vv{\cdot}\nabla)\bm\xi
+\Big[\frac{\phi_{\rm ref}^{\bm\xi}(J)}{J}\Big]_J'\vv{\cdot}\nabla J\,.
\end{align}\end{subequations}
 Making use of \eqref{DT-det} we find  
\begin{align}\nonumber
\int_\varOmega&\nabla p{\cdot}\vv\,\d\xx=
-\int_\varOmega\nabla[\phi_{\rm ref}^{\bm\xi}]_J'(J){\cdot}\vv\,\d\xx
=\int_\varOmega[\phi_{\rm ref}^{\bm\xi}]_J'(J){\rm div}\,\vv\,\d\xx
-\!\int_{\pl\varOmega}
[\phi_{\rm ref}^{\bm\xi}]_J'(J)\!\!\!\!\!\lineunder{\vv{\cdot}\nn}{$=0$}\!\!\!\!\d S
\\[-.5em]\nonumber
&\!\!\!\!\stackrel{\eqref{DT-det}}{=}
\int_\varOmega\frac{[\phi_{\rm ref}^{\bm\xi}]_J'(J)}J\DT J\,\d\xx
=\int_\varOmega\bigg(\Big[\frac{\phi_{\rm ref}^{\bm\xi}(J)}J\Big]_{\!J}'
+\frac{\phi_{\rm ref}^{\bm\xi}(J)}{J^2}\bigg)\DT J\,\d\xx
\\\nonumber
&\!\!\!\!\stackrel{\eqref{DT-det}}{=}
\int_\varOmega  
 \Big[\frac{\phi_{\rm ref}^{\bm\xi}(J)}J\Big]_{\!J}' \pdt J 
+\Big[\frac{\phi_{\rm ref}^{\bm\xi}(J)}J\Big]_{\!J}'\vv{\cdot}\nabla J
+\frac{\phi_{\rm ref}^{\bm\xi}(J)}J{\rm div}\,\vv\,\d\xx
\\\nonumber
&\!\!\!\!\stackrel{\eqref{calculus}}{=}\int_\varOmega
\pdt{}\Big(\frac{\phi_{\rm ref}^{\bm\xi}(J)}J\Big)
+\nabla\Big(\frac{\phi_{\rm ref}^{\bm\xi}(J)}J\Big){\cdot}\vv-\frac{\big[[\phi_{\rm ref}]_\XX'\big]^{\bm\xi}(J)\!}J{\cdot}
\Big(\!\!\!\!\!\lineunder{\pdt{\bm\xi}{+}(\vv{\cdot}\nabla)\bm\xi}{$=0$ due to \eqref{Euler3-fluid-selfgravit}}\!\!\!\!\!\Big)
+\frac{\phi_{\rm ref}^{\bm\xi}(J)}J{\rm div}\,\vv\,\d\xx
\\[-1.2em]\nonumber
&=\frac{\d}{\d t}\int_\varOmega\frac{\phi_{\rm ref}^{\bm\xi}(J)}J\,\d\xx
+\int_\varOmega{\rm div}\Big(\frac{\phi_{\rm ref}^{\bm\xi}(J)}J\,\vv\Big)
\,\d\xx
\\&=\frac{\d}{\d t}\int_\varOmega\frac{\phi_{\rm ref}^{\bm\xi}(J)}J\,\d\xx
+\int_{\pl\varOmega}
\frac{\phi_{\rm ref}^{\bm\xi}(J)}J\!\!\!\!\!\lineunder{\vv{\cdot}\nn}{$=0$}\!\!\!\!\,\d S
=\frac{\d}{\d t}\int_\varOmega\frac{\phi_{\rm ref}^{\bm\xi}(J)}J\,\d\xx\,.
\label{test-p}\end{align}
 Note that these computations require  sufficient  smoothness of
$\phi_{\rm ref}$  with respect to  $\XX$, so that
$[\phi_{\rm ref}]_\XX'$  is  integrable.   On the other hand,
from  the application point of view, it is desirable to consider 
sharp interfaces  between different materials,  which leads to jumps in
$\phi_{\rm ref}(\cdot,J):\Omega\to\R$.  Such material laws 
could also be rigorously  accounted for, at the expense of more refined
arguments,  cf.\ \cite{Roub22IFSV}  and Remark
\ref{rem:SharpInterfaces}.  

Altogether, under assumptions \eqref{stored-fluid} the energy balance now
reads
\begin{align}\nonumber
&\hspace*{0em}\frac{\d}{\d t}\Bigg(\int_\varOmega\!\!\!
\linesunder{\frac{\rho_{\rm ref}^{\bm\xi}}{2J}|\vv|^2}{kinetic}{energy}
\!\!\!\!+\!\!\!\!
  \morelinesunder{\frac{\phi_{\rm ref}^{\bm\xi}(J)}{J}}{actual}{stored}{energy}\!\!\!\!\!
+\!\!\!\!\!\!\!\!\!\!
\morelinesunder{\frac{\rho_{\rm ref}^{\bm\xi}}{J}\GRAVPOT}{\ \ energy of $\varrho$ in}{\ \ gravitational}{field $\GRAVPOT$}\!\!\!\!\!\!\!\!\!\d\xx
+\int_{\Universe}\!\!\!\!\!\!\!\!\morelinesunder{\frac{|\nabla\GRAVPOT|^2}{\TTT8\uppi\EEE\GRAVCONST}}{energy of}{gravitational}{field}\!\!\!\!\!\!\d\xx
+\int_{\Universe\setminus\varOmega}\!\!\!\!\!\!\!\!\morelinesunder{\varrho_{\rm ext}^{}\GRAVPOT}{energy of $\varrho_{\rm ext}^{}$}{\ \ \ in gravitational}{field $\GRAVPOT$}\!\!\!\!\!\!\!\!\!\!\!\!\!\d\xx\Bigg)
\\[-.0em]&\hspace{12em}
+\!\int_\varOmega
\!\!\!\!\linesunder{
\nu_1^{\bm\xi}|\EE(\vv)|^2
+\nu_2^{\bm\xi}|\Nabla\EE(\vv)|^{ q }_{_{_{}}}\!}{$=: \eta $ dissipation rate}{\ \ \ \ \ \ due to viscosity}
\!\!\!\d\xx
=\int_{\Universe\setminus\varOmega}\!\!\!\!\!\!\!\!\linesunder{\!\!\pdt{\varrho_{\rm
  ext}}\GRAVPOT}{power of}{external  forces }\!\!\!\!\!\!\!\d\xx\,.
\label{mech-engr-creep}
\end{align}

\begin{remark}[{\sl   Anisothermal  extension}]\label{rem-anisotherm1}\upshape
   In the spirit of Remark \ref{rem-anisotherm0}, also in case of
  assumption \eqref{stored-fluid} the model can be extended to the
   anisothermal  case by letting material  parameters  depend on the
  temperature $\theta$,  e.g., the  (referential) free energy  be 
  $\psi_{\rm ref}=\psi_{\rm ref}(\XX,J,\theta)$. Considering the split
  $\psi_{\rm ref}(\XX,J,\theta)=\phi_{\rm ref}(\XX,J)+\COUPLING_{\rm
    ref}(\XX,J,\theta)$ with $\gamma_{\rm ref}(\XX,J,0)=0$, the (actual)
  pressure $p$ in \eqref{Euler1-referential-selfgravit}  reads  
  $p=-[\phi_{\rm ref}^{\bm\xi}]_J'(\xx,J)-[\COUPLING_{\rm
    ref}^{\bm\xi}]_J'(\xx,J,\theta)$.  Hence, the ensuing heat equation is
\begin{align}\label{heat-eq+}
\pdt\W+{\rm div}\big(\vv\W{+}\jj\big)= \eta 
+ [\COUPLING_{\rm ref}^{\bm\xi}]_J'(J,\theta){\rm div}\,\vv
\ \ \ \text{ with }\ \ \W=\frac{\COUPLING_{\rm ref}^{\bm\xi}(  
J,\theta){-}\theta[\COUPLING_{\rm ref}^{\bm\xi}]_\theta'(J,\theta)}J \, ,
\end{align}
 where 
the heat production rate $ \eta  $ is  defined in 
\eqref{mech-engr-creep}.  This   anisothermal   version of the model again
reproduces the expected energetics and complies with the Clausius-Duhem
inequality. In particular, it can be checked to be thermodynamically
consistent. 
\end{remark}

\begin{remark}[{\sl State equation}]\label{rem-state-eq}\upshape
In fluid thermomechanics, the state equation relates density, pressure,
and temperature. Here, in view of Remark~\ref{rem-anisotherm1} and
 the fact that 
$\varrho=\rho_{\rm ref}/J$, this relation at a current material point reads as
$p=-\phi_{\rm ref}'(\rho_{\rm ref}/\varrho)
-[\COUPLING_{\rm ref}]_J'(\rho_{\rm ref}/\varrho,\theta)$. In the isothermal
situation, where $\COUPLING_{\rm ref}\equiv0$, the relation
$p=-\phi_{\rm ref}'(\rho_{\rm ref}/\varrho)$ represents a so-called
{\it isentropic state equation} while the fluid is  said to be  {\it barotropic},
i.e., its density depends only on pressure.
\end{remark}

\section{Analysis of  the viscoelastic fluid problem
}\label{sec-anal}

 We present here an existence result for weak solutions to system
\eqref{Euler-fluid-selfgravit}. Moreover, we prove that such weak
solutions fulfill the energy balance \eqref{mech-engr-creep}.
\DELETE{Note that in order to obtain  \eqref{mech-engr-creep},  as well as the
corresponding a-priori estimates, one needs a Poincar\'e inequality to
control the terms $\int_\varOmega\varrho\GRAVPOT\,\d\xx$,
$\int_{\Universe}\varrho_{\rm ext}\GRAVPOT\,\d\xx$, and
$\int_{\Universe}\pdt{}\varrho_{\rm ext}\GRAVPOT\,\d\xx$, which have
no sign, through $\int_{\R^3}|\nabla\GRAVPOT|^2\,\d\xx$. This 
once again asks $U$ to be bounded.}

We are interested in  an  initial-value problem for system 
\eqref{Euler-fluid-selfgravit}.  
The initial conditions  are prescribed as 
\begin{align}\label{IC-xi}
\vv|_{t=0}^{}=\vv_0\ \ \ \text{ and }\ \ \ 
\bm\xi|_{t=0}^{}=\bm\xi_0\,.
\end{align}

We will use the standard notation  for  Lebesgue and   Sobolev
spaces, namely $L^p(\varOmega;\R^n)$ for Lebesgue measurable functions
$\varOmega\to\R^n$ whose Euclidean norm  has integrable 
$p$-power, \linebreak
$W^{k,p}(\varOmega;\R^m)$ for functions from $L^p(\varOmega;\R^m)$ whose
 distributional derivatives  up to the order $k$ have their Euclidean norm integrable with
$p$-power, and $W_0^{k,p}(\varOmega;\R^m)$  for the  subspace of
$W^{k,p}(\varOmega;\R^m)$ of functions with zero  trace  on $\pl\varOmega$.
We also  use  $H^k=W^{k,2}$ and  $H_0^k=W_0^{k,2}$. The notation
$p^*$ will denote the  optimal exponent for  the embedding
$W^{1,p}(\varOmega)\subset L^{p^*}(\varOmega)$, i.e., $p^*=\TTT 3 \EEE
p/(\TTT 3 \EEE {-}p)$
for $p<\TTT 3 \EEE$ while $p^*\ge1$ arbitrary for $p=\TTT 3 \EEE$ or
$p^*=+\infty$ for $p>\TTT 3 \EEE$.
Moreover, for a Banach space
$X$ and for $I=[0,T]$, we will use the notation $L^p(I;X)$ for the Bochner
space of Bochner measurable functions $I\to X$ whose norm is in $L^p(I)$
while $W^{1,p}(I;X)$  indicates the  functions $I\to X$ whose distributional
derivative is in $L^p(I;X)$.  Moreover,  $C(\cdot)$ and $C^k(\cdot)$
will denote spaces of continuous and $k$-times continuously
differentiable functions,  respectively.  Eventually, $C_{\rm
  w}(I;X)$  and   ${\rm BV}(I;X)$
will denote the Banach space of weakly continuous functions $I\to X$
 and  functions with bounded variations, respectively, 
and $C_{\rm b}^{}(\cdot)$ will stand for  continuous  bounded functions.

 In the following, we use the symbol $C$ to indicate a generic
positive constant depending on data and possibly changing from line
to line. We impose the following assumptions on the data:
\COMMENT{Do we really want $\varOmega$ and $\Omega$ here??}
\begin{subequations}
 \label{ass}
  \begin{align}
    &\label{ass-omega}
    \Omega 
    \subset\R^3\ \text{\TTT is an \EEE open \TTT and \EEE bounded
      \TTT Lipschitz domain\EEE},
    \\\label{ass-phi}
    &\phi_{\rm ref}\in C^1( \ubarOmega  {\times}(0,+\infty)),\ \
    \exists\,\epsilon>0,\ \alpha>\TTT1\EEE
    \ \ \forall \XX\,{\in}\,\Omega,\ J>0:\ \
    \phi_{\rm ref}(\XX,J)\ge\frac\epsilon{J^{\alpha}}
    \,,
    \\&\label{ass-nu}
    \nu_1,\nu_2\in C( \ubarOmega )\,,
    \ \ \ 
    \min(\nu_1,\nu_2)>0\,,\ \ \
    q>3 \ \text{  (for $q$ occurring in \eqref{T+D})}
    \,,
    \\[-.0em]&\label{ass-rho}
    \rho_{\rm ref}\in C(\ubarOmega)\,,\ \ \ \ \rho_{\rm ref}>0\,,
    \ \ \ \GRAVCONST>0\,,\ \ \ 
    \\[-.0em]&\label{ass-rho-ext}
    \AAA \exists\, U\subset\R^3\text{ compact with } 
  \varOmega \cup {\rm sppt}\big(\varrho_{\rm ext}(t,\cdot)\big)\subset U, 
   \quad   \TTT  \text{ and }\ \ \varrho_{\rm ext}\in \AAA H^1(I;L^{2}(U))\,,\EEE
    \\[-.0em]&\label{ass-J0}
    \vv_0\in L^2(\varOmega;\R^3)\,,\ \ \ \bm\xi_0\in W^{2,r}(\varOmega;\R^3)\
    \text{ with }\ r>3\,,\ \ \ 1/\det({\nabla\bm\xi}_0)>0
    \ \text{ on }\barOmega\,.
\end{align}\end{subequations}
Note that $\phi_{\rm ref}(\XX,\cdot)$ need not be convex, hence
$\FF\mapsto\phi_{\rm ref}(\XX,\det\FF)$ need not be polyconvex.
\DELETE{In view of  the problem being independent of the actual value of
the constant $\GRAVPOT_\text{\sc b}^{}$, with no loss of generality we
choose $\GRAVPOT_\text{\sc b}^{}=0$ in eqref{BC-V} for simplicity.}

\begin{definition}[Weak solutions  of  system \eqref{Euler-fluid-selfgravit}]\label{def}
A  triplet  $(\vv,\bm\xi,\GRAVPOT)$ with 
$\vv\in C_{\rm w}(I;L^2(\varOmega;\R^3))$ $ \cap$ $ L^2(I;W^{2,q}(\varOmega;\R^3))$,
$\bm\xi\in C_{\rm w}(I;W^{2,r}(\varOmega;\R^3))
\cap W^{1,1}(I{\times}\varOmega;\R^3)$,
and $\GRAVPOT\in C_{\rm w}(I;H_0^{\TTT1\EEE}(\Universe))$
is a \emph{weak solution}  of  the boundary-value problem
\eqref{Euler-fluid-selfgravit} with the boundary conditions \eqref{BC} and
with the initial condition \eqref{IC-xi}
if
\begin{align}\nonumber
&\int_0^T\!\!\!\int_\varOmega\!\bigg(\big(\nu_1^{\bm\xi}\EE(\vv)
-\det(\nabla\bm\xi)\varrho_{\rm ref}^{\bm\xi}\vv{\otimes}\vv\big){:}\EE(\wt\vv)
+\nu_2^{\bm\xi}|\nabla\EE(\vv)|^{q-2}\nabla\EE(\vv)\Vdots\nabla\EE(\wt\vv)
-\det(\nabla\bm\xi)\varrho_{\rm ref}^{\bm\xi}\vv{\cdot}\pdt{\wt\vv}
\\[-.2em]&\hspace{4em}
+\bm[\phi_{\rm ref}^{\bm\xi}]_J'\Big(\frac1{\det(\nabla\bm\xi)}\Big){\rm div}\,\wt\vv
+\det(\nabla\bm\xi)\varrho_{\rm ref}^{\bm\xi}
\nabla\GRAVPOT{\cdot}\wt\vv\bigg)\,\d\xx \;\!\d t
=\int_\varOmega\frac{\varrho_{\rm ref}^{\bm\xi_0}\vv_0}{J_0}{\cdot}\wt\vv(0)\,\d\xx
\label{w-sln-u}
\end{align}
holds for all $\wt\vv\in C^\infty(I{\times}\barOmega;\R^3)$ with
$\wt\vv{\cdot}\nn=0$ and $\wt\vv(T)=\bm0$, 
\eqref{Euler3-fluid-selfgravit} hold a.e.\ on $I{\times}\varOmega$ with
$\bm\xi(0)=\bm\xi_0$ on $\varOmega$,
and \eqref{Euler2-fluid-selfgravit} holds
\TTT in the sense \eq{V-from-integral} \EEE a.e.\ on $I{\times}\Universe$.
\end{definition}

\begin{theorem}[Existence of solutions  of   \eqref{Euler-fluid-selfgravit}] 
\label{th:WeakSol}
Let  assumptions   \eqref{ass}  hold. 
Then:\\
\Item{(i)}{there exists a weak solution $(\vv,\bm\xi,\GRAVPOT)$
of system \eqref{Euler-fluid-selfgravit}.} 
\Item{(ii)}{Weak  solutions of system \eqref{Euler-fluid-selfgravit}  satisfy  the energy-dissipation balance
\eqref{mech-engr-creep}  when 
 integrated over  the  time interval $[0,t]$  for all  $t\in I$.
}
\end{theorem}

\def\EPS{\varepsilon}
\def\EPSk{{k}}

\begin{proof}
 The proof relies on a {\it semi-Galerkin} approximation
 and is divided into four steps.

\medskip\noindent{\it Step 1: semi-Galerkin approximation}.  
We perform a Galerkin approximation of the momentum equation
\eqref{Euler1-fluid-selfgravit}
for $\vv$  but leave  
the Poisson equation \eqref{Euler2-fluid-selfgravit} for $\GRAVPOT$ 
\TTT remaining continuous by relying on \eq{V-from-integral}\EEE.
\DELETE{the invertibility of the the Laplacian operator $-\Delta:H^2_0(\Universe)\to L^2(\Universe)$.}
 Similarly, we do not  discretize  
 the transport equation \eqref{Euler3-fluid-selfgravit} for $\bm\xi$ 
 and  rely on \cite[Lem.\,3.2]{RouSte??VSPS} for its weak solvability. 

 In order to approximate the momentum equation
\eqref{Euler1-fluid-selfgravit}, we introduce a family of 
nested finite-dimensional subspaces $\{\mathscr{V}_k\}_{k=0}^\infty$
whose union is dense in $W^{2,q}(\varOmega;\R^3)$.
Without loss of generality, we may assume $\vv_0\in\mathscr{V}_0$.

The global existence on the whole time interval \AAA $I=[0,T]$ \EEE of a
solution of such \DELETE{regularized and} a semi-discretized system, which will
be denoted by $(\vv_\EPSk,\bm\xi_\EPSk,\GRAVPOT_\EPSk)$, can be proved to exist
by the standard successive-prolongation argument, on the basis of the
uniform-in-time estimates proved below.

Setting  $J_\EPSk:=\det(\nabla\bm\xi_\EPSk)$ and
$\varrho_\EPSk:=\rho_{\rm ref}^{\bm\xi_\EPSk}/J_\EPSk$, we can rely on the
evolution-and-transport equations
\begin{align}\label{J-flow-k}
\pdt{J_\EPSk}=({\rm div}\,\vv_\EPSk)\,J_\EPSk-\vv_\EPSk{\cdot}\nabla J_\EPSk
\ \ \ \ \text{ and }\ \ \ \ \pdt{\varrho_\EPSk}=-{\rm div}(\varrho_\EPSk\vv_\EPSk)
\end{align}
with the initial conditions $J_\EPSk(0)=1/\det(\nabla\bm\xi_0)$ and
$\varrho_\EPSk(0):=\rho_{\rm ref}^{\bm\xi_0}/J_0$, respectively, cf.\ \eqref{DT-det}
and \eqref{cont-eq+}. Here,  we crucially used the fact that the
transport equation  
\eqref{Euler3-fluid-selfgravit}
is not discretized.

\medskip\noindent{\it Step 2: first a-priori estimates}.  
We test the Galerkin approximate versions of \eqref{Euler1-fluid-selfgravit}
and \eqref{Euler2-fluid-selfgravit} by $\vv_\EPSk$ and $\GRAVPOT_\EPSk$,
respectively, and use \eqref{J-flow-k}. The discretized velocity field
$\vv_\EPSk$ is in $L^{2}(I;W^{1,\infty}(\varOmega;\R^3))$, so that $J_\EPSk$,
which fulfills the non-discretized transport-and-evolution equation
\eqref{J-flow-k}, stays positive on $I{\times}\varOmega$; here assumption 
\eqref{ass-J0} is used, cf.\ \cite[Lem.\,3.2]{RouSte??VSPS}. By abbreviating
$\overline\nu_i=\min\nu_i(\XX)$ for $i=1,2$\TTT,
at each instant of time \EEE we find
\begin{align}\nonumber
&\hspace*{-1em}
\int_\varOmega\!\frac{\varrho_\EPSk}2|\vv_\EPSk|^2\!+\frac{ \epsilon}{\!J_\EPSk^{\alpha+1}(t)\!\!}
\,\d\xx
+\!\int_{\Universe}\!\!\frac{|\nabla\GRAVPOT_\EPSk(t)|^2\!\!}{\TTT8\uppi\EEE\GRAVCONST}\,\d\xx
+\!\int_0^t\!\!\int_\varOmega\!\overline\nu_1|\EE(\vv_\EPSk)|^2
+\overline\nu_2|\Nabla\EE(\vv_\EPSk)|^q
\,\d\xx \;\!\d t
\\&\nonumber\ \ \stackrel{\eqref{ass}}{\le}\int_\varOmega
\frac{\rho_{\rm ref}^{\bm\xi_\EPSk(t)}\!\!}{2J_\EPSk(t)}|\vv_\EPSk(t)|^2+
\frac{\phi_{\rm ref}^{\bm\xi_\EPSk(t)}(J_\EPSk(t))\!\!}{J_\EPSk(t)}\,\d\xx
+\!\int_0^t\!\!\int_\varOmega\!\nu_1^{\bm\xi_\EPSk}|\EE(\vv_\EPSk)|^2
+\nu_2^{\bm\xi_\EPSk}|\Nabla\EE(\vv_\EPSk)|^q
\,\d\xx \;\!\d t
\\&\ \ \stackrel{\eqref{mech-engr-creep}}{{\TTT\le\EEE}}
\int_0^t\!\!\int_{\TTT U\EEE}
\pdt{\varrho_{\rm ext}}\GRAVPOT_\EPSk\,\d\xx\,\d t
-\int_\varOmega
\frac{\rho_{\rm ref}^{\bm\xi_\EPSk(t)}}{J_\EPSk(t)}\GRAVPOT_\EPSk(t)\,\d\xx
-\!\int_{\TTT U\EEE}\!\varrho_{\rm ext}(t)V_\EPSk(t)\,\d\xx+\TTT C_0\EEE
\label{mech-engr-diff+}
\end{align}
with $\alpha$ from \eqref{ass-phi} \TTT and with $C_0$ depending on the
initial conditions qualified in \eq{ass-J0}.

For any time instant $t$ (omitted below for brevity), the gravitational
energy in \eq{mech-engr-diff+} can be estimated via \eqref{eq:V.lin.oper} as follows:
\begin{align} 
-\int_\varOmega\varrho_\EPSk\GRAVPOT_\EPSk\,\d\xx
\le \AAA \| \varrho_k\|_{L^2(\varOmega)} \| V_k\|_{L^2(\varOmega)} \leq C \|
\varrho_k\|_{L^2(\varOmega)} \|
\varrho_k{+} \varrho_{\rm ext}\|_{L^2(\varOmega)} \leq C+C\|
\varrho_k\|_{L^2(\varOmega)}^2  ,
\label{total-engr-selfgrav-accret-modif+}\end{align}
where $\varrho_\EPSk=\rho_{\rm ref}^{\bm\xi_\EPSk}/J_\EPSk$.
%
%
Using \eq{cont-eq+}, we continue to 
estimate \eq{total-engr-selfgrav-accret-modif+} as
$C+ C\,(\sup_{\varOmega}^{}\varrho_{\rm ref}^{})
\int_\varOmega1/J_k(\xx)^2\,\d\xx$, which eventually allows for  absorbing 
them in the left-hand side if $\phi(J)$ has a growth faster than $1/J$ for
$J\to0+$, as assumed in \eq{ass-phi}. The term
$\pdt{}\varrho_{\rm ext}\GRAVPOT_\EPSk$ in \eq{mech-engr-diff+} can be treated
by the Young and the Gronwall inequalities
if $\pdt{}\varrho_{\rm ext}\in L^2(I;L^2(\varOmega))$.
\EEE

From this,
we  obtain the a-priori bounds
\begin{subequations}\label{est-v-V-J}\begin{align}
\label{est-rv2}
&\big\|\sqrt{\varrho_\EPSk}\vv_\EPSk\big\|_{L^\infty(I;L^2(\varOmega;\R^3))}^{}\le C\,,
\\&\label{est-e(v)}
\|\EE(\vv_\EPSk)\|_{L^{\TTT q \EEE}(I; W^{1,q}(\varOmega;\R^{3\times 3}))}^{}\le C\,,
\\&\|\GRAVPOT_\EPSk\|_{L^\infty(I;H^1(\Universe))}^{}\le C
\,.
\intertext{By using \cite[Lem.\,3.2]{RouSte??VSPS}, from \eqref{J-flow-k}
for $J_\EPSk$ and assumption \eqref{ass-J0} we further obtain}
\label{est-J}
&\|J_\EPSk\|_{L^\infty(I;W^{1,r}(\varOmega))}^{}\le C\ \ \text{ and }\
\min J_\EPSk>1/C\,.
\intertext{From \eqref{Euler3-fluid-selfgravit}, i.e., $\pdt{}\bm\xi_\EPSk
  =-(\vv_\EPSk{\cdot}\nabla)\bm\xi_\EPSk$ and thus also
  $\pdt{}\nabla\bm\xi_\EPSk=-(\vv_\EPSk{\cdot}\nabla)\nabla\bm\xi_\EPSk
  -(\nabla\bm\xi_\EPSk)\nabla\vv_\EPSk$ and  taking advantage of the
                regularity  of the initial  value 
                $\nabla\bm\xi_0\in 
  W^{1,r}(\varOmega;\R^{\UUU 3\times 3\EEE})$,
   again  by \cite[Lem.\,3.2]{RouSte??VSPS} we also obtain}
\label{est-xi}
&\|\bm\xi_\EPSk\|_{L^\infty(I;W^{2,r}(\varOmega;\R^3))}^{}\le C\,.
\intertext{Since $\vv_\EPSk=(\sqrt{\varrho_\EPSk}\vv_\EPSk)\,
   (1/\sqrt{\varrho_\EPSk})  
  = (\sqrt{\varrho_\EPSk}\vv_\EPSk)\sqrt{J_\EPSk} / 
     \sqrt{\varrho_{\rm ref}^{\bm\xi_\EPSk}}$,
from \eqref{est-rv2} and \eqref{est-J} we eventually obtain}
&\label{est-v}
\|\vv_\EPSk\|_{L^\infty(I;L^2(\varOmega;\R^3))}^{}\le
\frac{\|J_\EPSk\|_{L^\infty(I\times\varOmega)}^{1/2}}{\min\varrho_{\rm ref}(\barOmega)^{1/2}}
\big\|\sqrt{\varrho_\EPSk}\vv_\EPSk\big\|_{L^\infty(I;L^2(\varOmega;\R^3))}^{}\le
  C\,. 
\intertext{\SSS Using \eq{mech-engr-diff+}, we can also bound 
$\nabla\EE(\vv_\EPSk)$ in $L^q(I{\times}\varOmega;\R^{3\times 3\times 3})$,
so that, by a generalized Poincar\'e inequality exploiting also \eq{est-v},
time-integrability in \eq{est-e(v)} is eventual improved to \EEE}
&\label{est-e(v)+}
\SSS\|\EE(\vv_\EPSk)\|_{L^q(I; W^{1,q}(\varOmega;\R^{3\times 3}))}^{}\le C\,.\EEE
\end{align}\end{subequations}

From  \eqref{J-flow-k} for $J_k$ 
we obtain a bound for $\pdt{}J_k$ in $L^q(I;L^r(\varOmega))$ and,
similarly, from $\pdt{}\bm\xi_k=-(\vv_k{\cdot}\nabla)\bm\xi_k$, we obtain a
bound for $\pdt{}\bm\xi_k$ in $L^q(I;L^r(\varOmega;\R^3))$.  This
additionally provides a bound on  
$\pdt{}\nabla\bm\xi_k$,  which we will not use, however.
From $\pdt{}\GRAVPOT_k\TTT(\xx)=
\int_{U}\pdt{}(\varrho_k(\wt\xx){+}\varrho_{\rm ext}(\wt\xx))/|\xx{-}\wt\xx|\,\d\wt\xx$ \TTT
we can see that $\nabla\pdt{}\GRAVPOT_k(\xx)=
\int_{U}\pdt{}(\varrho_k(\wt\xx){+}\varrho_{\rm ext}(\wt\xx))
(\xx{-}\wt\xx)/|\xx{-}\wt\xx|^3\,\d\wt\xx$. \EEE
With $\pdt{}\varrho_k\in L^q(I;L^r(\varOmega))$ and
$\pdt{}\varrho_{\rm ext}\in \AAA L^2(I;
\TTT L^2(U)\EEE)$, we
obtain a bound for $\pdt{} V_k$ in $ \AAA L^2 \EEE (I;\TTT H^2(U)\EEE)$ \AAA by
using $q>3$ from \eqref{ass-nu} and \eqref{eq:V.lin.oper}. \EEE

\medskip\noindent{\it Step 3: limit passage with $k\to\infty$}.
By the Banach selection principle, we select a weakly* convergent
subsequence and $(\varrho,\vv,\bm\xi,J,\GRAVPOT)$ such that
\begin{subequations}\label{Euler-weak}\begin{align}\label{Euler-weak-rho}
&\varrho_\EPSk\to\varrho&&\text{weakly* in }\
L^\infty(I;W^{1,r}(\varOmega))\,\cap\,W^{1,q}(I;L^r(\varOmega))\,,
\\&\vv_\EPSk\to\vv&&\text{weakly* in }\
L^\infty(I;L^{\SSS2\EEE}(\varOmega;\R^3))\,\cap\,
L^{\SSS q\EEE}(I;W^{2,q}(\varOmega;\R^3))\,,
\\&\label{Euler-weak-xi}
\bm\xi_\EPSk\to \bm\xi&&\text{weakly* in }\ L^\infty(I;W^{2,r}(\varOmega;\R^3))\,\cap\,
W^{1,q}(I;L^{\SSS r\EEE}(\varOmega;\R^3))\,,
\\&J_\EPSk\to J&&\text{weakly* in }\ L^\infty(I;W^{1,r}(\varOmega))\,\cap\,
W^{1, q }(I;L^r(\varOmega))\,,
\\ \label{eq:Vk-V}&\GRAVPOT_\EPSk\to\GRAVPOT&&\text{weakly in }\
\AAA H^1(I; H^2(U)) ,\EEE
\\ \label{eq:nabla.Vk}
  &\hspace*{-0.81em}\AAA \nabla \GRAVPOT_\EPSk\to\nabla \GRAVPOT&&
 \AAA \text{weakly in }\ H^1(I; H^1(\R^3;\R^3) )\,.
\intertext{\AAA For the last two relations we use the linear relation
  $V=\mathsf V[\varrho{+}\varrho_{\rm ext}]$ and \eqref{eq:V.lin.oper}. \EEE
Recalling that $r>3$, by  the Aubin-Lions Lemma, the convergences 
(\ref{Euler-weak}a,c,d) are also strong:}
&J_k\to J&&\text{strongly in$\ C(I{\times}\barOmega)\,$},
\\&\varrho_k\to\varrho=\rho_{\rm ref}^{\bm\xi}/J &&\text{strongly in
$\ C(I{\times}\barOmega)\,$, and}\\
&\bm\xi_\EPSk\to\bm\xi&&\text{strongly in
$\ C(I{\times}\barOmega;\R^3)$}.
\intertext{Moreover, $\varrho_\EPSk\vv_\EPSk\to\varrho\vv$ and $\vv_\EPSk\to\vv$
strongly in $L^c(I{\times}\varOmega;\R^3)$ for all $1\le c<4$, cf.\
\cite{RouSte??VSPS} for details. From the mentioned strong convergence of
$\varrho_\EPSk$ and \eqref{eq:V.lin.oper}, we
further obtain}
&\nabla\GRAVPOT_\EPSk\to \nabla\GRAVPOT
&&\text{strongly in }\ C(I;L^2(\R^3;\R^3))\,.
\label{nabla-V-strong}
\end{align}
\end{subequations}

We now use the Galerkin approximation of the momentum equation
\eqref{Euler1-fluid-selfgravit} tested by
$\wt\vv=\vv_\EPSk-\wt\vv_k$ where
$\wt\vv_k:I\to\mathscr{V}_k$  is  an approximation of $\vv$ such that
$\wt\vv_k\to\vv$ strongly in $L^\infty(I;L^2(\varOmega;\R^{\TTT3\EEE}))$
and $\Nabla\EE(\wt\vv_k)\to\Nabla\EE(\vv)$ strongly in
$L^q(I{\times}\varOmega;\R^{\TTT 3\times 3\times 3\EEE})$ for
$k\to\infty$.
\TTT 
Using the continuity equation from \eqref{J-flow-k} one
directly checks that
\begin{equation}\label{eq:useme}
  \varrho_{\EPSk}\DT\vv_{\EPSk}= \frac{\partial}{\partial t} (\varrho_{\EPSk}
\vv_{\EPSk})+{\rm div} (\varrho_{\EPSk}\vv_{\EPSk}{\otimes}
\vv_{\EPSk}).
\end{equation}
Integrating on $I{\times}\varOmega$ equation \eqref{eq:referlater} written for $\varrho_\EPSk$ and
$\vv_\EPSk$ and using \eqref{eq:useme} we obtain
\begin{align*}
  &\int_\varOmega \frac{\varrho_\EPSk(T)}2 |\vv_\EPSk(T)|^2\,\d\xx -
\int_\varOmega\!\frac{\varrho_0}2|\vv_0|^2\,\d\xx = \int_0^T \frac{\d}{\d t}
\int_\varOmega 
    \left(\frac{\varrho_\EPSk}2|\vv_{\EPSk}|^2\right)\,\d\xx\,\d t\\
  &\qquad\stackrel{\eqref{eq:referlater}}{=} \int_0^T \!\!\!\int_\varOmega
\varrho_\EPSk \DT\vv_{\EPSk}{\cdot}\vv_{\EPSk}\,\d\xx\,\d t
\stackrel{\eqref{eq:useme}}{=}\int_0^T \!\!\!\int_\varOmega \Big(\frac{\partial}{\partial t}  (\varrho_{\EPSk}
\vv_{\EPSk})+{\rm div} (\varrho_{\EPSk}\vv_{\EPSk}{\otimes}
\vv_{\EPSk})\Big){\cdot}\vv_{\EPSk}\,\d\xx\,\d t.
\end{align*}
The latter implies that \EEE
\begin{align}\nonumber
\!\!\int_\varOmega\!\frac{\varrho_\EPSk(T)}2\big|\vv_\EPSk(T){-}\vv(T)\big|^2\d\xx
&=\int_0^T\!\!\!\int_\varOmega\Big(\pdt{}(\varrho_\EPSk\vv_\EPSk)
+{\rm div}(\varrho_\EPSk\vv_\EPSk{\otimes}\vv_\EPSk)
\Big){\cdot}\vv_\EPSk\,\d\xx\!\;\d t
\\[-.1em]&
\quad +\!\int_\varOmega\!\frac{\varrho_0}2|\vv_0|^2\!-\varrho_\EPSk(T)\vv_\EPSk(T){\cdot}\vv(T)
+\frac{\!\varrho_\EPSk(T)}2|\vv(T)|^2\d \xx\,.
\label{Euler-one-substitution}
\end{align}
\SSS Letting $c_q>0$ be such that 
$c_{q}|H-\widetilde H|^q\le(|H|^{q-2}H-|\widetilde H|^{q-2}\widetilde
H)\Vdots(H-\widetilde H)$ holds for all $H \in {\mathbb R}^{3\times
  3\times 3}$, we can estimate \EEE
\begin{align}\nonumber
&\frac{\min\rho_{\rm ref}^{}(\overline\Omega)}{2\|J_\EPSk(T)\|_{L^\infty(\varOmega)}^{}}
\big\|\vv_\EPSk(T){-}\vv(T)\big\|_{L^2(\varOmega;\R^{ 3})}^2\!\!
                         \TTT
                         +\bar\nu_1\|\EE(\vv_\EPSk{-}\vv)
                         \|_{L^2(I{\times}\varOmega;\R^{3\times 3})}^2
                         \EEE
                         +\bar\nu_2\,c_{q}\,\|\nabla\EE(\vv_\EPSk{-}\vv)
\|_{L^q(I{\times}\varOmega;\R^{3\times 3\times 3})}^q
\\&\nonumber
\hspace*{1em}\overset{\text{convexity}}\leq
 \int_\varOmega\frac{\varrho_\EPSk(T)}2\big|\vv_\EPSk(T){-}\vv(T)\big|^2\,\d\xx
+\int_0^T\!\!\!\int_\varOmega\!\bigg(
\nu_1^{\bm\xi_\EPSk}\ee(\vv_\EPSk{-}\SSS\vv\EEE){:}\ee(\vv_\EPSk{-}\SSS\vv\EEE)
  \\[-.2em]&\hspace*{10em}\nonumber
 +\nu_2^{\bm\xi_\EPSk}\big(|\nabla\EE(\vv_\EPSk)|^{q-2}\nabla\EE(\vv_\EPSk)
-|\nabla\EE(\vv)|^{q-2}\nabla\EE(\vv)\big)\Vdots
  \nabla\EE(\vv_\EPSk{-}\vv)\bigg)\,\d\xx \;\!\d t
 \\[-.5em]&\hspace*{2em}=\nonumber
 \int_0^T\!\!\!\int_\varOmega\bigg(p_\EPSk{\rm div}(\vv_\EPSk{-}\wt\vv_k)
-\varrho_\EPSk\nabla\GRAVPOT_\EPSk{\cdot}(\vv_\EPSk{-}\wt\vv_k)
-\nu_1^{\bm\xi_\EPSk}\EE(
 \SSS\vv\EEE){:}\EE(\vv_\EPSk{-} \SSS \vv \EEE )
\\[.2em]&\nonumber\hspace{3em}
{}-\nu_2^{\bm\xi_\EPSk}|\nabla\EE(
\SSS\vv\EEE)|^{ q-2 }\nabla\EE(
\SSS\vv\EEE)\Vdots\nabla\EE(\vv_\EPSk{-} \SSS \vv\EEE)
 +\Big(\pdt{}(\varrho_\EPSk\vv_\EPSk)
+{\rm div}(\varrho_\EPSk\vv_\EPSk{\otimes}\vv_\EPSk)\Big)
{\cdot}\wt\vv_k\bigg)\,\d\xx \;\!\d t
\\[-.1em]&\hspace{3em}
 +\int_\varOmega\bigg(\frac{\varrho_0}2|\vv_0|^2
 -\varrho_\EPSk(T)\vv_\EPSk(T){\cdot}
 \SSS\vv\EEE(T) +\frac{\varrho_\EPSk(T)}2|
\SSS\vv\EEE(T)|^2\bigg)\,\d\xx+\mathscr{R}_k \,
\label{strong-hyper+}\end{align}
with $p_\EPSk=-[\phi_{\rm ref}^{\bm\xi_\EPSk}]_J'(J_\EPSk)$ \SSS coming from
using \eqref{Euler-one-substitution}, 
inserting the Galerkin approximation of \eqref{Euler1-fluid-selfgravit} tested
by $\vv_k{-}\wt\vv_k$, and defining the remainder $ \mathscr{R}_k$ as \EEE
\begin{align}
 \nonumber
\mathscr{R}_k& \SSS =
\int_0^T\!\!\!\int_\varOmega\!
\nu_1^{\bm\xi_\EPSk}
 \EE(\vv_\EPSk)  {:}\EE(\wt\vv_k{-}\vv)
 +\nu_2^{\bm\xi_\EPSk}|\nabla\EE(\vv_\EPSk)|^{ q-2 }\nabla\EE(\vv_\EPSk)\Vdots
  \nabla\EE(\wt\vv_k{-}\vv)\,\d\xx \;\!\d t
\end{align}
\TTT The aim is now to show that the right-hand side in \eqref{strong-hyper+}
converges to $0$. The first four term converge to $0$ by strong-weak converges,
while the fifth term has a limit that cancels with the last integral by using
\eqref{Euler-one-substitution} once again, see \cite{RouSte??VSPS, Roub22TVSE}
for details. \EEE
The remainder term $\mathscr{R}_k$ in
\eqref{strong-hyper+} converges to $0$ due to the strong approximation
properties of the 
$\wt\vv_k$ for $\vv$. 

Thus, we obtain the strong convergence
\begin{subequations}\label{strong-conv}\begin{align}
&&&\vv_\EPSk\to
\SSS\vv\EEE&&\text{strongly in $L^q(I;W^{2,q}(\varOmega;\R^{ 3}))$}
\intertext{ together with $\vv_\EPSk(T)\to
\SSS\vv\EEE(T)$ in $L^2(\varOmega;\R^{\TTT3\EEE})$.
In fact, by performing this computation at a generic time $t\in (0,T)$
instead of $T$, we obtain}
&&&\vv_\EPSk(t)\to\vv(t)\!\!\!\!&&\text{strongly in $\,L^2(\varOmega;\R^{ 3})\,$ for any $t\in I$.}&&&&&&
\end{align}\end{subequations}

Owing to these convergences, the passage to the limit 
in the semi-discrete system is  straightforward. In particular,  
the limit is a weak solution  of the system in the sense of 
Definition~\ref{def}. \DELETE{By differentiating
\eqref{Euler2-fluid-selfgravit} in time and taking into account  
that $\pdt{}\varrho\in L^p(I;L^r(\varOmega))$, from  
the assumption $\pdt{}\varrho_{\rm ext}\in L^1(I;L^{6/5}(\varOmega))$
we have that  $\GRAVPOT\in W^{1,1}(I;
W^{2,6/5}_0(\Universe))$ and, in particular, $\GRAVPOT\in C_{\rm w}(I;
W^{2,6/5}_0(\Universe))$.
Here, we used the  classical elliptic $W^{2,p}$-regularity
theory
as well as the fact that the
right-hand side $\varrho_{\rm ext}\in W^{1,1}(I;L^{6/5}(\Universe))$
is fixed, so that the limit time derivative is integrable.}

\medskip\noindent{\it Step 4: energy-dissipation balance}.  To
conclude the proof, we now check that  
the calculations leading to \eqref{mech-engr-selfgravit}
in the variant \eqref{mech-engr-creep}, in particular \eqref{test-p},
are indeed legitimate. Here, we simply refer  to  \cite[Sec.\,3]{Roub22IFSV}
or \cite[Sec.\,3]{RouSte??VSPS}  where this check has been already
performed. With respect to these references, additional care has to be
given here in order to obtain \eqref{calculus-selfgravit}, since  
the mechanical load $\nabla p+\varrho\nabla\GRAVPOT$ in the momentum equation
\eqref{Euler1-fluid-selfgravit} has to be shown to be in duality with
$\vv\in L^\infty(I;L^2(\varOmega;\R^3))\cap L^q(I;W^{2,q}(\varOmega;\R^3))$.
This follows 
\TTT from the fact that \EEE $p$ is in $L^\infty(I;W^{1,r}(\varOmega))$
so that $\nabla p$ is \TTT surely \EEE in $L^1(I;L^2(\varOmega;\R^3))$, and
the restriction of $\nabla\GRAVPOT$ to $\varOmega$ is \TTT surely \EEE in
$L^\infty(I;L^2(\varOmega;\R^3))$.
\DELETE{Note in addition that $\pdt{}\GRAVPOT$ needs to be in
$L^1(I;L^6(\Universe))$ in order to be in duality with the equation
\eqref{Euler2-fluid-selfgravit}. This however follows from  
$\pdt{}\GRAVPOT\in L^1(I;W^{2,6/5}_0(\Universe))$  which has been 
proved in Step~3.}
\TTT By a comparison in \eqref{Euler2-fluid-selfgravit},
$\Delta\GRAVPOT\in H^1(I;L^2_{\rm loc}(\Universe))$ and vanishes outside $U$ and
is in duality with $\pdt{}\GRAVPOT\in L^2(I;H^1(U))$. Together with
$\nabla\GRAVPOT\in H^1(I;H^1(\Universe;\R^3))$,
this makes the formal calculus \eq{calculus-selfgravit} analytically
legitimate. \EEE
\end{proof}

\begin{remark}[{\sl Long time  scales}]\label{rem-long-time}  
\upshape
Having in mind the self-gravitational differentiation of planets during
long time scales, when the initial configuration is successively forgotten, 
\TTT the discussion of \EEE the validity of the above estimates for
$T\to\infty$ is relevant. The constants in the bounds in
(\ref{est-v-V-J}d-\SSS{}g\EEE) depend on the regularity of the initial
conditions and are  possibly (exponentially) increasing in time.  
On the other hand, estimates (\ref{est-v-V-J}a,c) are controlled
by the material properties and  can pass to the limit $T\to\infty$
(at least after neglecting \TTT possible \EEE effects of the moving external
mass $\varrho_{\rm ext}$). Moreover, from \eqref{mech-engr-diff+} we have  
an estimate of $1/J$ in $L^\infty( 0,+\infty;L^{1+\alpha}(\varOmega))$ \TTT
so that \EEE also the linear momentum
$\varrho\vv=\sqrt{\varrho}\vv(\rho_{\rm ref}^{\bm\xi}/J)^{1/2}$ is bounded.
Specifically for $\alpha\TTT>1\EEE$ and using \eqref{est-rv2}, it holds
\begin{align*}
&\|\varrho\vv\|_{L^\infty( 0,+\infty;L^{\TTT4/3\EEE}(\varOmega;\R^3))}^{}
\le C\TTT\EEE
\big\|\sqrt{\varrho}\vv\big\|_{L^\infty( 0,+\infty ;L^2(\varOmega;\R^3))}^{}
\Big\|\frac1{\sqrt J}\Big\|_{L^\infty( 0,+\infty ;L^{\TTT4\EEE}(\varOmega))}.
\end{align*}
\end{remark}

\begin{remark}[{\sl Pressure dependent viscosities}]\upshape
Often, viscosity coefficients depend (beside temperature, not considered
here) also on the pressure \cite{RaSaVe09OBAF} and may  vary, 
in particular during phase transitions in some materials
\cite{SSSY77VDPE,SNBS12IPDV}. 
This can be  taken into account  by  letting  $\nu_1=\nu_1(\XX,J)$
and  $\nu_2=\nu_2(\XX,J)$  depend on $J$. 
\end{remark}

\begin{remark}[{\sl Sharp interfaces}]\upshape \label{rem:SharpInterfaces}
In many applications (and in particular in planetary geophysics),
the solid is composed by very different materials. Correspondingly 
the reference data $\phi_{\rm ref}^{}(\cdot,J)$ and $\rho_{\rm ref}^{}$ as well
as the viscosities $\nu_1(\cdot,J)$ and $\nu_2(\cdot,J)$ are  ideally
discontinuous with respect to $\XX$, as opposed to 
\eqref{ass}. This makes the substitutions with $\bm\xi$ analytically more
complicated, \TTT because \EEE one cannot use the continuity of the implied 
Nemytski\u{\i} mapping (composition), 
as actually used in the above proof. Instead, one  has  to modify the
free-slip boundary condition \eqref{BC}  and assume  a stick condition
$\vv=\bm0$ on $\pl\varOmega$.  This, together with the local invertibility
$\det(\nabla\bm\xi(t))>0$, would ensure the global invertibility of
$\bm\xi$, eventually allowing to apply a change of variable. The reader is
referred to \cite[Sec.\,4]{Roub22IFSV} for additional details in a similar
setting. Let us just record that this approach hinges on a bound for the
distorsion, which in turn ensures that referential interfaces between the
regions occupied by different materials remain of null measure in the
actual deformed configuration.
\end{remark}

\begin{remark}[{\sl\TTT The surroundings $U$\EEE}]\upshape\label{rem:U}
\TTT The bounded set $U$ in \eq{ass-rho-ext} selects a ``neighbourhood of
influence''  for
the self-graviting planet and thus the objects in the whole Universe which
are responsible for the considered tidal effects caused by $\varrho_{\rm ext}$.
In fact, it also influences the values of the gravitational potential
$\GRAVPOT$ which, in reality, is determined rather vaguely. Let us \EEE
note that \SSS it is \EEE well-known that the actual value of \TTT this \EEE
potential on the Earth surface depends on whether one considers only
a single planet, the whole Solar System, or the whole galaxy Milky
Way,  and \TTT its value is then \EEE roughly 60 MJ/kg, 900 MJ/kg, or more than
130 GJ/kg, respectively.
\end{remark}

\section{Multicomponent materials}\label{sec-multi}

The different regions of planets and moons are actually composed by many
distinct  materials (cf., e.g., \cite{Cond16EEPS,ShTuOl04MCEP}). 
These materials,  may undergo pressure-dependent
chemical reactions, combined with diffusion.  In this section, we
extend the model by including the description of the different
constituent of the solid by means of  
a concentration vector $\cc=(c_1,...,c_n)$. 
 We  assume that the viscoelastic  response of the solid
 depends on the composition,  namely, we let
$\nu_1=\nu_1(\XX,J,\cc)$, $\nu_2=\nu_2(\XX,J,\cc)$, as well as  $\phi_{\rm ref}=\phi_{\rm ref}(\XX,J,\cc)$.  On the other hand, we assume the
mass density to be independent of $\cc$, 
 so that $\varrho$ is still  determined by \eqref{rho=rho0/detF}, cf.\ Remark~\ref{rem=mass} below.
Without loss of generality, we take  the number 
$n$ of constituents to be  the same  in all regions  (hence $n$ is independent of $\XX$). 
Of course, the components of $\cc$ are to
be non-negative and  to  satisfy $\sum^n_{i=1}c_i=1$ a.e.\
 in  $I{\times}\varOmega$.  In other words, 
$\cc$  takes values 
in the so-called Gibbs' simplex
$$\triangle^+_1:=\{(c_1,...,c_n)\in\R^n;\ \ \sum^n_{i=1}c_i=1\ \
 \text{and}  \ \   \forall i:\ c_i\ge0\}.$$

The  single-component  system \eqref{Euler-fluid-selfgravit} is
 then  expanded  to its multi-component variant  as
\begin{subequations}\label{Euler-fluid-selfgravit+}\begin{align}
\nonumber
&\varrho\DT\vv=
{
\rm div}\big(
\nu_1^{\bm\xi}(\cc)\EE(\vv)-{\rm div}\big(\nu_2^{\bm\xi}(\cc)|\nabla\EE(\vv)|^{q-2}\nabla\EE(\vv)\big)\big)
-\nabla p
          -     \varrho\nabla\GRAVPOT\ \ \ 
 \\& \hspace{7.5em}\text{ where }\ \
p= -\frac{[\phi_{\rm ref}^{\bm\xi}]_J'(J,\cc)}{J} 
\ \ \text{ and }\ \ \varrho=\frac{\rho_{\rm ref}^{\bm\xi}}J
\ \text{ with }\ J=\frac1{\det(\nabla\bm\xi)}
\,,
\label{Euler1-fluid-selfgravit+}
\\\label{Euler3-diff}
&\DT{\cc}={\rm div}\big(\bbM^{\bm\xi}(J,\cc)
\nabla\bm\mu\big) -\rr^{\bm\xi}(J,\cc)
\\[-.7em]\label{Euler3-diff+}
&\qquad\qquad\qquad\qquad
\ \ \ \ \ \text{ with }\ \
\bm\mu\in\frac{[\phi_{\rm ref}^{\bm\xi}]_{\cc}'(J,\cc)}{J}+\bm{N}(\cc)
\,,
\\[-.5em]\label{Euler3-fluid-selfgravit+}
&\DT{\bm\xi}=\bm0,\hspace{7.2em}\text{on }\ \varOmega\,, 
\\\label{Euler2-fluid-selfgravit+}
&\Delta\GRAVPOT=\TTT4\uppi\EEE\GRAVCONST\big(
\varrho{+}\varrho_{\rm ext}^{}\big)\qquad\
\text{ on }\ \Universe\,,
\end{align}\end{subequations}
where $\bm{N}
(\cc)$ in \eqref{Euler3-diff+} denotes the normal cone to the convex set
$\triangle^+_1$
at $\cc$.  In particular, $\bm{N}
(\cdot):\R^n\rightrightarrows\R^n$ is  a maximal monotone 
set-valued mapping.  In fact, the term 
$\bm{N}(\cc)$ in \eqref{Euler3-diff+}  contributes a Lagrange 
multiplier (= a ``pressure'')  corresponding  to the  constraints
$\sum^n_{i=1}c_i=1$ and $c_i \geq 0$.  Such multiplier
ensures the validity of the  constraints throughout the  
evolution.  In this context, the use of such a multiplier dates at least
back to 
E and Palffy-Muhoray \cite{EPal97PSIS},  who nonetheless 
generalized an (essentially) 1D-model  by  De\,Gennes
\cite{DeGe80DFSD}.  We also refer to \cite{OttE97TDIF} for a discussion of
local versus a nonlocal mixture models, where our approach corresponds to the
nonlocal model with multiplier $\lambda(t,\xx)\bm1\in \bm{N} (\cc)$ 
with $\bm1=(1,\ldots,1)\in \R^n$ associated with
the constraint $\cc(t,\xx) \in \triangle^+_1$, see also Remark \ref{rem-M}.

 Relation \eqref{Euler3-diff} features  the $n {\times} n$ mobility  matrix
$\bbM=\bbM(\XX,J,\cc)$ and the chemical-reaction rate $\rr=\rr(\XX,J,\cc)$.
The mobility matrix $\bbM(\XX,J,\cc)$ is  assumed  to be
symmetric  in order  to comply with the
Onsager principle.  In addition, it is assumed to be  positive
semi-definite to comply with the  Clausius-Duhem  inequality and
thus the 2nd-law of thermodynamics. For analytical reasons,  we
assume $\bbM$ to be 
uniformly positive definite,  as this allows  to control the chemical-potential
gradient and, indirectly, also the concentration gradients.
The mass conservation within chemical reactions  imposes that the
reaction rates 
 $\rr=(r_1,\dots,r_n)$  satisfy 
the condition
\begin{align}
\label{eq:rr.cdot1}
\forall(\XX,J,\cc)\in \Omega{\times}\R^+{\times}\triangle^+_1:\ \
&\sum_{i=1}^n r_i(\XX,J,\cc)=0
\,.
\end{align}
 Note that we are following here the  
phenomenological approach by Eckart  and  Prigogine
\cite{Ecka40TIP,Prig47ETPI} by assuming that all component have the
same velocity $\vv$. A 
 less  phenomenological,  truly rational-thermodynamical  alternative
would be to assume that each constituent has its own velocity, as in
the Truesdell \cite{True68BTM} approach.  

The boundary conditions \eqref{BC} are  complemented by the
boundary condition $\nn{\cdot} \bbM^{\bm\xi}(J,\cc)\nabla\bm\mu=\bm0$ on 
$\pl\varOmega$ for \eqref{Euler3-diff}, expressing that
there is no flux of the constituents across $\pl\varOmega$. 
The energetics  of the model can be deduced  
as in Section~\ref{sec-fluids}, now
combined with \eqref{Euler3-diff} tested by $\bm\mu$, which gives 
\begin{align}
\int_\varOmega\DT\cc{\cdot}\bm\mu\,\d\xx
= -  \int_\varOmega\rr^{\bm\xi}(J,\cc){\cdot}\bm\mu
+\nabla\bm\mu{:}\bbM^{\bm\xi}(J,\cc)\nabla\bm\mu\,\d\xx\,,
\end{align}
and further using
\eqref{Euler3-diff+} tested by $\DT\cc$. 
We modify \eqref{test-p} to be merged with a part of \eqref{Euler3-diff+}
tested by $\DT\cc$. Specifically,  also using  \eqref{transport-xi} similarly
as in \eqref{test-p}, we have
\begin{align}\nonumber
\int_\varOmega\bm\mu{\cdot}\DT\cc&+\nabla p{\cdot}\vv
\,\d\xx
=\int_\varOmega \bm\mu{\cdot}\DT\cc-p{\cdot}{\rm div}\vv
=\int_\varOmega\!\frac{[\phi_{\rm ref}^{\bm\xi}]_\cc'(J,\cc)}J{\cdot}\DT\cc
+\frac{[\phi_{\rm ref}^{\bm\xi}]_J'(J,\cc)}J\DT J\,\d\xx
\\\nonumber
&=\int_\varOmega\bigg(\Big[\frac{\phi_{\rm ref}^{\bm\xi}(J,\cc)}J\Big]_{\!J}'
+\frac{\phi_{\rm ref}^{\bm\xi}(J,\cc)}{J^2}\bigg)\DT J
+\frac{[\phi_{\rm ref}^{\bm\xi}]_\cc'(J,\cc)}J{\cdot}\DT\cc\,\d\xx
\\\nonumber
&=\int_\varOmega\pdt{}\Big(\frac{\phi_{\rm ref}^{\bm\xi}(J,\cc)}J\Big)
+\Big[\frac{\phi_{\rm ref}^{\bm\xi}(J,\cc)}J\Big]_{\!J}'\vv{\cdot}\nabla J
+\Big[\frac{\phi_{\rm ref}^{\bm\xi}(J,\cc)}J\Big]_{\!\cc}'{\cdot}(\vv{\cdot}\nabla)\cc
+\frac{\phi_{\rm ref}^{\bm\xi}(J,\cc)}J{\rm div}\,\vv
\\\nonumber
&=\int_\varOmega
\pdt{}\Big(\frac{\phi_{\rm ref}^{\bm\xi}(J,\cc)}J\Big)
+\nabla\Big(\frac{\phi_{\rm ref}^{\bm\xi}(J,\cc)}J\Big){\cdot}\vv
+\frac{\phi_{\rm ref}^{\bm\xi}(J,\cc)}J{\rm div}\,\vv\,\d\xx
\\
&=\frac{\d}{\d t}\int_\varOmega\!\frac{\phi_{\rm ref}^{\bm\xi}(J,\cc)}J\,\d\xx
+\int_\varOmega\!{\rm div}\Big(\frac{\phi_{\rm ref}^{\bm\xi}(J,\cc)}J\,\vv\Big)\,\d\xx
=\frac{\d}{\d t}\int_\varOmega\!\frac{\phi_{\rm ref}^{\bm\xi}(J,\cc)}J\,\d\xx\,.
\label{test-p-c}\end{align}
The  remaining  term arising from $\bm{N}(\cc)$ in \eqref{Euler3-diff+}
tested by $\DT\cc$  can be proved to vanish, provided that
$\cc(0,\xx)\in\triangle_1^+$.  Indeed, let $\delta_{\triangle_1^+}$ be the
indicator function of $\triangle_1^+$ and $\bm{\eta}\in \NN(\cc)$ a.e.
in
$I{\times}\varOmega$.  By comparison from \eqref{Euler3-diff+}  one
finds that
$\bm{\eta}\in L^2(I{\times}\varOmega;\R^n)$ (under the assumption
\eqref{ass-phi+} below). Hence, the classical chain rule for convex
functions (see, for instance \cite[Prop.\,XI.4.11]{Visi96MPT})
ensures 
\begin{align}\label{chain}
\bm{\eta}{\cdot}\DT\cc=\bm{\eta}{\cdot}\pdt{\cc} +
   \bm{\eta}{\cdot}(\vv{\cdot}\nabla)\cc 
= \DT{\overline{\delta_{\triangle_1^+}(\cc)}}=0
\end{align}
a.e.\ in $I{\times}\varOmega$, the last equality following from the
fact that $\cc \in \triangle_1^+$ a.e. in $I{\times}\varOmega$. 
In fact, by comparison from \eqref{Euler3-diff+} one can see that
$\bm{\eta}\in L^2(I{\times}\varOmega;\R^n)$ under the assumption
\eqref{ass-phi+} below, as needed in \cite{Visi96MPT} for a rigorous
proof of \eqref{chain}.  

We hence deduce the energy-dissipation balance
\begin{align}\nonumber
&\hspace*{-.5em}\frac{\d}{\d t}\Bigg(\int_\varOmega\!\!\!
\linesunder{\frac{\rho_{\rm ref}^{\bm\xi}}{2J}|\vv|^2}{kinetic}{energy}
\!\!\!\!+\!\!\!\!
  \morelinesunder{\frac{\phi_{\rm ref}^{\bm\xi}(J,\cc)}{J}}{actual}{stored}{energy}\!\!\!\!\!
+\!\!\!\!\!\!\!\!\!\!
\morelinesunder{\frac{\rho_{\rm ref}^{\bm\xi}}{J}\GRAVPOT}{\ \ energy of $\varrho$ in}{\ \ gravitational}{field $\GRAVPOT$}\!\!\!\!\!\!\!\!\!\d\xx
+\int_{\Universe}\!\!\!\!\!\!\!\!\morelinesunder{\frac{|\nabla\GRAVPOT|^2}{\TTT8\uppi\EEE\GRAVCONST}}{energy of}{gravitational}{field}\!\!\!\!\!\!\d\xx
+\int_{\Universe\setminus\varOmega}\!\!\!\!\!\!\!\!\!\!\!\!\!\morelinesunder{\varrho_{\rm ext}^{}\GRAVPOT_{_{_{{_{}}}}}}{energy of $\varrho_{\rm ext}^{^{^{}}}$}{\ \ \ in gravitational}{field $\GRAVPOT$}\!\!\!\!\!\!\!\!\!\!\!\!\!\d\xx\Bigg)
\\[-.0em]&\hspace{.5em}\nonumber
+\!\int_\varOmega
\!\!\!\!\linesunder{
\nu_1^{\bm\xi}(\cc)|\EE(\vv)|^2\!
+\nu_2^{\bm\xi}(\cc)|\Nabla\EE(\vv)|^{ q}_{_{_{}}}\!}{dissipation rate}{due to viscosity}
\!\!\!+\!\!\!\!\linesunder{\bbM^{\bm\xi}(J,\cc)\nabla\bm\mu{:}\nabla\bm\mu}{
dissipation rate}{due to diffusion}
\!\!\!\!+\!\!\!\!\!\!\!\morelinesunder{\rr^{\bm\xi}(J,\cc){\cdot}\bm\mu}{dissipation
  rate}{ due to reactions}\!\!\!\!\!\! \d\xx
=\int_{\Universe\setminus\varOmega}\hspace{-3em}\morelinesunder{\!\!\pdt{\varrho_{\rm ext}}
\GRAVPOT}{ gravitational power\ \ }{of external\ \ }{mass}\hspace{-3em}\d\xx\,.
\\[-1.5em]\label{mech-engr-selfgravit+}
\end{align}

The initial conditions \eqref{IC-xi}  now include a prescription for
$\cc$ and read  
\begin{align}\label{IC-c}
\vv|_{t=0}^{}=\vv_0\,,\ \ \ \ \ \bm\xi|_{t=0}^{}=\bm\xi_0\,,
\ \ \ \text{ and }\ \ \ \cc|_{t=0}^{}=\cc_0\in\triangle^+_1\,.
\end{align}

\begin{definition}[Weak solutions of the system \eqref{Euler-fluid-selfgravit+}]
\label{def2}
The quintuple $(\vv,\cc,\bm\mu,\bm\xi,\GRAVPOT)$ is a \emph{weak
  solution}  of  the
boundary-value problem for  the system  \eqref{Euler-fluid-selfgravit+} with the boundary
conditions \eqref{BC} together with $(\nn{\cdot}\nabla)\bm\mu=\bm0$ on
$\pl\varOmega$ and with the initial conditions \eqref{IC-c}
if $(\vv,\bm\xi,\GRAVPOT)$ is as in Definition~\ref{def} with
\eqref{w-sln-u} holding 
with
$\nu_1^{\bm\xi}=\nu_1^{\bm\xi}(\cc)$ and
$\nu_2^{\bm\xi}=\nu_2^{\bm\xi}(\cc)$
and if $\cc\in L^2(I;H^1(\varOmega;\R^n))$ and
$\bm\mu\in L^2(I;H^1(\varOmega;\R^n))$
satisfy  $0\le c_i\le 1$ and $\sum_{i=1}^nc_i=1$ a.e.\ on $I{\times}\varOmega$,
the integral identity
\begin{align}
&\int_0^T\!\!\!\int_\varOmega\bbM^{ \bm \xi}( J,\cc)\nabla\bm\mu{:}\nabla\widetilde{\bm\mu}
+\big((\vv{\cdot}\nabla)\cc  +  \rr^{\bm\xi}(J,\cc)\big){\cdot}\widetilde{\bm\mu}
-\cc{\cdot}\pdt{\widetilde{\bm\mu}}\,\d\xx \;\!\d t=\int_\varOmega\cc_0{\cdot}
\widetilde{\bm\mu}(0)\,\d\xx
\label{w-sln-c}
\end{align}
 holds  for any $\widetilde{\bm\mu}\in H^1(I{\times}\varOmega;\R^n)$ with
 $\widetilde{\bm\mu}(T)=\bm0$,  and the inclusion \eqref{Euler3-diff+} holds a.e.\ on
$I{\times}\varOmega$.
\end{definition}

 In order to treat the multi-component case, the assumptions on data have to
 be specified as follows,  where we use the short-hand \TTT notation \EEE
$\mathfrak{D}:=\ubarOmega{\times}(0,+\infty){\times}\R^n$:  \COMMENT{THE NOTATION IN \eq{ass-m} WAS CHANGED TO AVOID ``G'' WHICH HAS A DIFFERENT USAGE}
\begin{subequations}\label{ass+}\begin{align}\nonumber
&
\phi_{\rm ref}(\XX,J,\cc)=\phi_{\rm r,0}(\XX,J)+\phi_{\rm r,1}(\XX,J,\cc),\ \
\phi_{\rm r,0}
\in C^1(\ubarOmega{\times}(0,+\infty))\,,\ \ 
\\&\nonumber\qquad
\phi_{\rm r,1}\in C^1(\mathfrak{D} )
\ \ \text{ with }\  [\phi_{\rm r,1}]_\cc'  \in C^1(\mathfrak{D}) 
 \ \  \text{and} \ \  [\phi_{\rm r,1}]''_{\cc\cc} \in
 C_{\rm b}(\mathfrak{D};\R^{n\times n}), 
\\&\nonumber\qquad
\exists\,\epsilon>0,\ \alpha>
\TTT1\EEE\,,\ \ \forall \XX\,{\in}\,\Omega,\ J>0,\ \cc\in\R^n:\ \
\phi_{\rm r,0}(\XX,J)\ge\frac\epsilon{J^{\alpha}}\,,
\\[-.4em]&\qquad\phi_{\rm r,1}(\XX,J,\cc)\ge0\,,\ \ \
\big|[\phi_{\rm r,1}]_J'(\XX,J,\cc)\big|\le\frac1{\!\epsilon J^{(\alpha+1)/2}}\,,\, 
\text{ and }\ 
[\phi_{\rm r,1}]_{\cc\cc}''(\XX,J,\cc)\ge
\epsilon
\,,
\label{ass-phi+}
\\& 
\bbM\in C_{\rm b}^{}(\mathfrak{D};\R_{\rm sym}^{(n\times3)^2})\,,\ 
\ 
\label{ass-m}
    \quad\forall \TTT A\EEE\in {\mathbb R}^{n\times 3}:\ \  
   \inf\nolimits_{(\XX,J,\cc)\in\mathfrak{D}}^{}
\TTT A\EEE^\top\!{:}\bbM(\XX,J,\cc){:}\TTT A\EEE\ge\epsilon|\TTT A\EEE|^2\,,
\\
&\nonumber  
\rr(\XX,J,c)= \bbK(\XX,J,\cc) [\phi_{\rm r,1}]'_\cc(\XX,J,\cc) \ \text{ with  some  }\ 
\bbK\in C_{\rm b}^{}(\mathfrak{D};\R_{\rm sym}^{n\times n}), \ \ \bbK\bm1
=\bm0, \text{ and }
                                  \\[-.3em]&\hspace*{6em}  
                                             \forall\, { \bm\mu}\in\R^n,\ 
       \sum_{i=1}^n\mu_i=0:\ \ 
     \inf\nolimits_{(\XX,J,\cc)\in\mathfrak{D}}^{}{ \bm\mu}{\cdot}\bbK(\XX,J,\cc){
                                             \bm\mu}\ge\epsilon |\bm\mu|^2,
\label{ass-r}
\\
&\label{ass-nu+}
\nu_1,\nu_2\in C(\ubarOmega{\times}\R^n)\,,
\ \ \ 
\mbox{$\inf_{\ubarOmega{\times}\R^n}^{}$}\min(\nu_1,\nu_2)>0\,,
\ \ \ q>3\,.
\end{align}
\end{subequations}
Note that \eqref{ass-phi+} in particular  implies   that
$\phi_{\rm ref}(\XX,J,\cdot)$  is  uniformly convex  with
respect to $\cc$. For the purposes of the mathematical study, this  implies that
the diffusion equation defined via \eqref{Euler3-diff}-\eqref{Euler3-diff+} is
parabolic. The fact that the reaction terms $\rr$ can be written in the form 
$\bbK \, [\phi_{\rm r,1}]'_\cc$ (under the condition of detailed balance) was
observed in \cite{Miel11GSRD}, see also \cite{MaaMie20MCRS} and Remark~\ref{rem:Stoich}  below.  The  condition $\bbK \bm1=\bm0$  in
\eqref{ass-r}  guarantees \eqref{eq:rr.cdot1}. 
 Note also that assumptions are  formulated 
also for $\cc\in\R^n\setminus\triangle^+_1$, which   
is required as we  will  use an 
 exterior-penalty technique below. 

\begin{theorem}[Existence of  weak solutions  of  system \eqref{Euler-fluid-selfgravit+}]\label{thm-1}
Let  assumptions   \eqref{ass-omega}, \eqref{ass-rho},
\eqref{ass-J0},  
and \eqref{ass+}  hold.  Then:\\
\Item{(i)}{there exists a weak solution  of  system \eqref{Euler-fluid-selfgravit+}.}
\Item{(ii)}{
  Weak  solutions of system \eqref{Euler-fluid-selfgravit+}  
  satisfy the energy-dissipation balance
\eqref{mech-engr-selfgravit+} 
when   integrated over
 the  time interval $[0,t]$ with any $t\in I$.}
\end{theorem}

\def\EPS{\varepsilon}
\def\EPSk{{\varepsilon k}}
\def\LAM{\lambda}

\begin{proof}
   We argue by approximation and subdivide the proof into 
  seven  subsequent  steps.  
In addition to  a  semi-Galerkin approximation  in the
spirit of 
Section~\ref{sec-anal},
 now used   also for \eqref{Euler3-diff}, we perform
a regularization of the
indicator function $\delta_{\triangle_1^+}(\cdot)$  by an 
exterior  penalization.

\medskip\noindent{\it Step 1: regularized problem}.  
The multivalued mapping  $\bm{N}
(\cdot)$ in \eqref{Euler3-diff+} is regularized by  introducing  a penalization 
$\mathcal{P}_\EPS $ of the indicator function of the 
Gibbs' simplex $\triangle_1^+$ defined here,
for any $\EPS>0$, by
\begin{align*}
&\mathcal{P}_\EPS(\cc)=\frac1{2\EPS}\sum_{i=1}^n\min(0,c_i)^2+
\frac1{2\EPS}\Big(\sum_{i=1}^nc_i-1\Big)^2\,.
\end{align*}
Note that $\mathcal{P}_\EPS$ is convex  and  continuously
differentiable.  As $\EPS \to 0$ one has that $\mathcal{P}_\EPS\to
\delta_{\triangle_1^+}$ pointwise and increasing. 
 
 The regularized problem is obtained by replacing  $\bm{N}$ in
\eqref{Euler3-diff+} by the derivative
$\mathcal{P}_\EPS'$.  
The inclusion \eqref{Euler3-diff+} turns  into the   equation
\begin{align}
\label{Euler3-diff++}
\bm\mu=
\frac{[\phi_{\rm ref}^{\bm\xi}]_{\cc}'(J,\cc)}{J}+\mathcal{P}_\EPS'(\cc)\,.
\end{align}

 This choice also regularizes the reaction terms $\rr$  as  
\begin{equation}
  \label{eq:rEPS.1}
  \rr^{\bm\xi}_\EPS = \bbK^{\bm\xi} \bm\mu = 
\bbK^{\bm\xi}\Big(\frac1J [\phi_{\rm ref}^{\bm\xi}]_{\cc}'(J,\cc) + \mathcal{P}_\EPS'(\cc)\Big).
\end{equation}

The weak solution of  the  boundary-value problem for the system
(\ref{Euler-fluid-selfgravit+}a,b,d,e) and \eqref{Euler3-diff++} with the boundary
conditions \eqref{BC} together with $(\nn{\cdot}\nabla)\bm\mu=\bm0$ on
$\pl\varOmega$ and the initial condition \eqref{IC-c} will be denoted
by $(\vv_\EPS,\cc_\EPS,\bm\mu_\EPS,\bm\xi_\EPS,\GRAVPOT_\EPS)$. 
Its existence is proved in Step 4 below. 

\medskip\noindent{\it Step 2: semi-Galerkin approximation}.  
We perform a 
Galerkin approximation of  the momentum equation \eqref{Euler-fluid-selfgravit+}
for $\vv$ as in Section~\ref{sec-anal}
and of the diffusion equation \eqref{Euler3-diff} for $\cc$. 
On the other hand, we do not approximate in space
 the transport equation \eqref{Euler3-fluid-selfgravit+} for $\bm\xi$  
 but rather rely on  \cite[Lem.\,3.2]{RouSte??VSPS} for its weak 
 solvability. Moreover, 
we do not approximate  neither  the Poisson equation \eqref{Euler2-fluid-selfgravit+}
 nor   the regularized nonlinear equation \eqref{Euler3-diff++}.

Specifically, we again use a nested finite-dimensional subspaces
$\{\mathscr{V}_k\}_{k=0}^\infty$
for the momentum equation \eqref{Euler1-fluid-selfgravit+}.
 For the 
Galerkin approximation of 
the  diffusion  equation \eqref{Euler3-diff}  we use  
a second collection of  nested
finite-dimensional subspaces $\{\mathscr{Z}_k\}_{k=0}^\infty$ whose union
is dense in $H^1(\varOmega;\R^n)$.
Without loss of generality, we may assume $\vv_0\in\mathscr{V}_0$
and $\cc_0\in\mathscr{Z}_0$.

 We directly substitute $\bm\mu$ from \eqref{Euler3-diff++} into
 \eqref{Euler3-diff} to obtain a parabolic equation for $\cc_\EPS$.
It should be emphasized that the physically motivated tests leading to the
discrete energy-dissipation balance \eqref{mech-engr-selfgravit+}
cannot be  performed at  the Galerkin-discretization level because 
$\DT\cc_\EPSk=\pdt{}\cc_\EPSk+(\vv_\EPSk{\cdot}\nabla)\cc_\EPSk$ is not 
a legitimate  test for  the  discretized equation for $\bm\mu_\EPSk$.
 This calls for implementing another estimation strategy.  Moving from
 \eqref{Euler3-diff++} we have that 
\begin{align}
\!\!\nabla\bm\mu&=\Big(\frac{[\phi_{\rm ref}^{\bm\xi}]_{\cc\cc}''(J,\cc)\!\!}{J}
+\mathcal{P}_\EPS''(\cc)\Big)\nabla\cc
+\Big(\frac{[\phi_{\rm ref}^{\bm\xi}]_{J\cc}''(J,\cc)\!}{J}
-\frac{[\phi_{\rm ref}^{\bm\xi}]_{\cc}'(J,\cc)}{J^2}\Big)\nabla J
+\frac{[[\phi_{\rm ref}]_{\XX\cc}'']^{\bm\xi}(J,\cc)\!}{J}\,\nabla\bm\xi
\label{nabla-mu}
\end{align}
so that,  by substituting \eqref{Euler3-diff++} into  
\eqref{Euler3-diff} 
 we obtain the  semilinear parabolic equation for
$\cc_\EPSk$:
\begin{align}\nonumber
&\pdt{\cc_\EPSk}-{\rm div}\big(M_\EPS^{\bm\xi_\EPSk}(J_\EPSk,\cc_\EPSk)\nabla\cc_\EPSk\big)
={\rm div}\Big(R^{\bm\xi_\EPSk}_{\EPS}(J_\EPSk,\cc_\EPSk)\nabla J_\EPSk
\\[-.4em]&\nonumber\hspace{17em}
{+}S^{\bm\xi_\EPSk}_{\EPS}(J_\EPSk,\cc_\EPSk)\nabla\bm\xi_\EPSk\Big)
-(\vv_\EPSk{\cdot}\nabla)\cc_\EPSk
  -  \rr^{\bm\xi_\EPSk}(J_\EPSk,\cc_\EPSk)
\\[.2em]&\nonumber\hspace{5.3em}\text{with }\ R_{\EPS}(\XX,J,\cc):= M_\EPS(\XX, J, \cc)
\Big(\frac{[\phi_{\rm ref}]_{J\cc}''(\XX,J,\cc)}{J}
-\frac{[\phi_{\rm ref}]_{\cc}'(\XX,J,\cc)}{J^2}\Big)\,,\ \ \ 
\\&\nonumber\hspace{8em}S_{\EPS}(\XX,J,\cc):= M_\EPS(\XX, J, \cc)
\frac{[\phi_{\rm ref}]_{\XX\cc}''(\XX,J,\cc)}{J}\,,\ \ \ \text{and }\ 
\\&\hspace{8em}
M_\EPS(\XX, J, \cc):=\bbM(\XX,J,\cc)\Big(
\frac{[\phi_{\rm ref}]_{\cc\cc}''(\XX,J,\cc)}{J}+\mathcal{P}_\EPS''(\cc)\Big)^{-1}.
\label{diff-c}
\end{align}

The  global existence on the whole time interval $I=[0,T]$ of a solution
of such regularized and semi-discretized system, which  we denote 
by 
$$
(\vv_\EPSk,\cc_\EPSk,\bm\xi_\EPSk,\GRAVPOT_\EPSk):I\to \mathscr{V}_k\times
\mathscr{Z}_k\times W^{2,r}(\varOmega;\R^{\TTT3\EEE})\times H^1(\Universe),
$$
results from a standard successive-prolongation argument, on the basis of
the uniform-in-time estimates proved below.

\medskip\noindent{\it Step 3: first a-priori estimates}. 
 We first test 
(\ref{Euler-fluid-selfgravit+}a,e) separately. 
Specifically, we 
test
the Galerkin discretization
of \eqref{Euler1-fluid-selfgravit+} by $\vv_\EPSk$
 and test 
\eqref{Euler2-fluid-selfgravit+} by $\GRAVPOT_\EPSk$, also 
using  \eqref{Euler3-fluid-selfgravit+}. The discretized velocity field
$\vv_\EPSk$ is in $L^{2}(I;W^{1,\infty}(\varOmega;\R^{\TTT3\EEE}))$ so that
$J_\EPSk=1/\!\det(\nabla{\bm\xi_\EPSk})$,
which fulfills the non-discretized transport-and-evolution
equation \eqref{J-flow-k}, stays positive on $I{\times}\varOmega$. 
Here, \cite[Lem.\,3.2]{RouSte??VSPS} has been used, \TTT exploiting
the \EEE assumption \eqref{ass-J0}. 

By arguing as in \eqref{mech-engr-diff+}, abbreviating
$\overline\nu_i=\inf\nu_i(\XX,\cc)$ with $i=1,2$, and     
using also \eqref{test-p} for $\phi_{\rm r,0}$ in place of $\phi_{\rm ref}$
 and  assumption \eqref{ass-phi+},  we obtain  
\begin{align}\nonumber
&\hspace*{0em}
\int_\varOmega\!\frac{\varrho_\EPSk}2|\vv_\EPSk|^2\!+\frac\epsilon{\!J_\EPSk^{\alpha+1}(t)\!\!}
\,\d\xx
+\!\int_{\Universe}\frac{|\nabla\GRAVPOT_\EPSk(t)|^2\!\!}{\TTT8\uppi\EEE\GRAVCONST}\,\d\xx
+\!\int_0^t\!\!\int_\varOmega\!\overline\nu_1|\EE(\vv_\EPSk)|^2
+\overline\nu_2|\Nabla\EE(\vv_\EPSk)|^q
\,\d\xx \;\!\d t
\\&\nonumber\ \stackrel{\eqref{ass}}{\le}\!\!\int_\varOmega
\frac{\rho_{\rm ref}^{\bm\xi_\EPSk(t)}\!\!}{2J_\EPSk(t)}|\vv_\EPSk(t)|^2+
\frac{\phi_{\rm r,0}^{\bm\xi_\EPSk(t)}(J_\EPSk(t))\!\!}{J_\EPSk(t)}\,\d\xx
\\&\qquad\qquad\qquad\ \ \nonumber
+\!\int_{\Universe}\frac{|\nabla\GRAVPOT_\EPSk(t)|^2}{\TTT8\uppi\EEE\GRAVCONST}\,\d\xx
+\!\int_0^t\!\!\int_\varOmega\!\nu_1^{\bm\xi_\EPSk}(\cc_\EPSk)|\EE(\vv_\EPSk)|^2
+\nu_2^{\bm\xi_\EPSk}(\cc_\EPSk)|\Nabla\EE(\vv_\EPSk)|^q
\,\d\xx \;\!\d t
\\&\nonumber\ \stackrel{\eqref{mech-engr-diff+}}{\le}\!\!
\int_0^t\!\bigg(\!\int_{U}\!\pdt{\varrho_{\rm ext}}\GRAVPOT_\EPSk\,\d\xx
+\int_{\varOmega}
\big[\phi_{\rm r,1}^{\bm\xi}(J_\EPSk,\cc_\EPSk)\big]_J'{\rm div}\,\vv_\EPSk
\,\d\xx\bigg)\,\d t
\\[-.4em]&\hspace{16em}
-\!\!\int_\varOmega
\frac{\rho_{\rm ref}^{\bm\xi_\EPSk(t)}}{J_\EPSk(t)}\GRAVPOT_\EPSk(t)\,\d\xx
-\!\int_{U}\!\!\varrho_{\rm ext}(t)V_\EPSk(t)\,\d\xx\TTT+C_0\,.\EEE
\end{align}
From this, \DELETE{taking $\delta>0$ sufficiently small} using \TTT
\eq{total-engr-selfgrav-accret-modif+} \EEE and the Gronwall inequality, we
obtain the a-priori bounds
\begin{subequations}\label{est-v-V-J+}\begin{align}
&\|\vv_\EPSk\|_{L^q(I;W^{2,q}(\varOmega;\R^3))}^{}\le C_\EPS, 
\\&\|\GRAVPOT_\EPSk\|_{L^\infty(I;H^1(\TTT U\EEE))}^{}\le C_\EPS\,.
\intertext{From the equations for $\bm\xi$ and
$J$  and  
assumption \eqref{ass-J0},  by  using
\cite[Lem.\,3.2]{RouSte??VSPS} we also obtain}
&\|\bm\xi_\EPSk\|_{L^\infty(I;W^{2,r}(\varOmega;\R^3))}^{}\le C_\EPS\,,\ \ \ \ 
\|J_\EPSk\|_{L^\infty(I;W^{1,r}(\varOmega))}^{}\le C_\EPS\,,\ \text{ and }\
\min J_\EPSk>1/C_\EPS\,.
\label{est-xi-eps-k}\end{align}\end{subequations}
It should be  noted 
that $C_\EPS$ here depends possibly on $\EPS$ but  is independent of 
$k$. From \eqref{est-xi-eps-k}, we also  obtain a 
bound for $\varrho_\EPSk=\rho_{\rm ref}^{\bm\xi_\EPSk}/J_\EPSk$ in
$L^\infty(I;W^{1,r}(\varOmega))$.

Next, we test the parabolic equation \eqref{diff-c} by $\cc_\EPSk$.
From \eqref{est-v-V-J+}, we have that  $R^{\bm\xi_\EPSk}$  is bounded in
$L^\infty (I{\times}\varOmega;  
\R^{n\times3}))$. Using the Green formula for
$$\int_\varOmega((\vv{\cdot}\nabla)\cc){\cdot}\cc\,\d\xx
=-\frac12\int_\varOmega({\rm div}\,\vv)|\cc|^2\d\xx$$  and the assumption
\eqref{ass-r}, this test gives
\begin{align}\nonumber
&
\int_\varOmega\!\frac{|\cc_\EPSk(t)|^2}2\,\d\xx
+\!\int_0^t\!\!\int_\varOmega M_\EPS^{\bm\xi_\EPSk}(J_\EPSk,\cc_\EPSk)\nabla\cc_\EPSk{:}\nabla\cc_\EPSk\,\d\xx \;\!\d t
\\\nonumber
&\quad=
\!\int_0^t\!\!\int_\varOmega
R^{\bm\xi_\EPSk}_{\EPS}(J_\EPSk,\cc_\EPSk) \nabla J_\EPSk {:}\nabla\cc_\EPSk
- \rr_\EPS^{\bm\xi_\EPSk} (J_\EPSk,\cc_\EPSk){\cdot}\cc_\EPSk
-\frac12({\rm div}\,\vv_\EPSk)|\cc_\EPSk|^2\!\,\d\xx \;\!\d t
\\&\nonumber\quad\le\!\int_0^t\!\bigg(\delta\|\nabla\cc_\EPSk\|_{L^2(\varOmega;\R^{n\times3})}^2
+\frac1{4\delta}\|R^{\bm\xi_\EPSk}_{\EPS}(J_\EPSk,\cc_\EPSk) 
 \nabla J_\EPSk \|_{L^2(\varOmega;\R^{n\times3})}^2
\\[-.4em]&\qquad\qquad
+ C_\EPS \big(|\varOmega|+\|\cc_\EPSk\|_{L^2(\varOmega;\R^n)}^2\big)
+\|{\rm div}\,\vv_\EPSk\|_{L^\infty(\varOmega)}^{}
\|\cc_\EPSk\|_{L^2(\varOmega;\R^n)}^2\bigg)\,\d t\,.
\end{align}
 Here we estimated $|\rr^{\bm\xi}_\EPS|= \big|\bbK^{\bm\xi}\,\big( \frac1{J_\EPS}
[\phi_{\rm r,1}]'_\cc + \mathcal P_\EPS'(\cc)\big) \big| \leq C_\EPS(1{+}|\cc|)$ $  $ by
using that $|\bbK|\leq C$ from \eqref{ass-r}, $|[\phi_{\rm
r,1} ]'_\cc|\leq C(1{+}|\cc|)$ from \eqref{ass-phi+}, and $|\mathcal{P}'_\EPS(\cc)| \leq
C(1{+}|\cc|)/\EPS$ and by exploiting estimate \eqref{est-xi-eps-k}.  

By the uniform positive definiteness of $M_\EPS(\cdot,\cdot)$, choosing
$\delta>0$ sufficiently small and exploiting the Gronwall inequality, we
also obtain
\begin{align}
\|\cc_\EPSk\|_{L^\infty(I;L^2(\varOmega;\R^n))\,\cap\,L^2(I;H^1(\varOmega;\R^n))}^{}\le C\,.
\end{align}
From \eqref{Euler3-diff++} and \eqref{nabla-mu}, we obtain the bound for
$\bm\mu_\EPSk$ in $L^\infty(I;L^2(\varOmega;\R^n))
\,\cap\,L^2(I;H^1(\varOmega;\R^n))$   uniformly in $k$, depending possibly
on $\varepsilon>0$. 

\medskip\noindent{\it Step 4: limit passage  for  $k\to\infty$}.
By the Banach selection principle, we select a weakly* convergent
subsequence and
$(\varrho_\EPS,\vv_\EPS,\bm\xi_\EPS,J_\EPS,\GRAVPOT_\EPS,\cc_\EPS,\bm\mu_\EPS)$
such that the convergences \eqref{Euler-weak} (still for fixed value
$\EPS>0$) hold, together with 
\begin{subequations}\begin{align}
&\cc_\EPSk\to\cc_\EPS&&\text{weakly* in }\ L^\infty(I;L^2(\varOmega;\R^n))
\,\cap\,L^2(I;H^1(\varOmega;\R^n))\,,
\\&\bm\mu_\EPSk\to\bm\mu_\EPS&&\text{weakly* in }\ L^\infty(I;L^2(\varOmega;\R^n))
\,\cap\,L^2(I;H^1(\varOmega;\R^n))\,.
\end{align}\end{subequations}
From \eqref{diff-c}, we  deduce a bound on 
$\pdt{}\cc_\EPSk$
in the respective semi-norms induced by the Faedo-Galerkin discretization
by the finite-dimensional subspaces $\mathscr{Z}_k$. By the (generalized)
Aubin-Lions theorem, we  obtain the 
strong convergence $\cc_\EPSk\to\cc_\EPS$ in
$L^{s}(I{\times}\varOmega;\R^n)$ for any $1\le s<10/3$. From
\eqref{Euler3-diff++}, 
such strong convergence also  holds  for
$\bm\mu_\EPSk\to\bm\mu_\EPS$.

 Adapting the argument leading to \eqref{Euler-weak}--\eqref{strong-conv} 
towards the weak formulation of \eqref{Euler3-diff} with
\eqref{Euler3-diff++} is then easy.

\medskip\noindent{\it Step 5: ``physically motivated'' a-priori estimates}.  
 As $\cc_\EPS$  is not discretized, we can test \eqref{Euler3-diff++}
by $\DT\cc_\EPS=\pdt{}\cc_\EPS+(\vv_\EPS{\cdot}\nabla)\cc_\EPS$. Together with
the other tests, we thus obtain  the analogous  energy-dissipation balance
to \eqref{mech-engr-selfgravit+}, now for the $\EPS$-solution and with the
additional left-hand-side term
$\frac{\d}{\d t}\int_\varOmega\mathcal{P}_\EPS(\cc_\EPS)\,\d\xx$. This delivers 
a-priori estimates, similarly as those obtained in \eqref{mech-engr-diff+},
with the additional energy estimate 
\begin{equation}
  \label{eq:calPeps}
  \int_\varOmega \mathcal{P}_\EPS\big(\cc_\EPS(t,\xx)\big)\,\d\xx \leq C \quad
  \text{for all } t \in I.
\end{equation}
 Let us introduce now the  orthogonal projection $\bbP:\R^n \to
\R^n$ to the subspace
$\{\bm\mu\in \R^n\,|\, \sum_1^n \mu_i=0\}$, namely,
$$(\bbP \bm\mu)_i = \mu_i - \frac1n \sum_{j=1}^n \mu_j = \frac1n
\sum_{j=1}^n (\mu_i-\mu_j) \quad \text{for} \ i=1,\dots,j,$$
so that
$$ \bbP\bm\mu = \bm\mu - \Big( \frac1n \sum_{j=1}^n \mu_j
\Big){\bm 1} = \bm \mu -\frac1n (\bm \mu {\cdot} \bm 1) \bm 1 .$$
Taking into account \eqref{ass-r} so that specifically $\bbK{\bm1}={\bm0}$,
the reaction term in the estimate gives  
\begin{align*} 
&\int_0^T\!\!\! \int_{\varOmega} \rr_\EPS^{\bm\xi_\EPS} 
(J_\EPS,\cc_\EPS){\cdot}\bm\mu_\EPS\,\d\xx \,\d t  = \int_0^T\!\!\!
  \int_\varOmega  \bm\mu_\EPS {\cdot} \bbK\bm\mu_\EPS \,\d\xx\d t
   = \int_0^T\!\!\!
  \int_\varOmega  \bbP\bm\mu_\EPS {\cdot}\bbK\,\bbP\bm\mu_\EPS \,\d\xx\d t
\geq  \epsilon \int_0^T \!\!\!\int_{\varOmega} |\bbP\bm\mu_\EPS|^2
  \,\d \xx \,\d t.
\end{align*}
We thus obtain estimates  as in   
\eqref{est-v-V-J+}, but for $\EPS$-solution, together with additional estimates
\begin{subequations}
  \label{est-c-eps}
 \begin{align}
   \label{est-mu-c-eps}
& \|\bbP\bm\mu_\EPS\|_{L^2(I{\times}\varOmega;\R^n)}^{}+  \|\nabla\bm\mu_\EPS\|_{L^2(I{\times}\varOmega;\R^{n\times3})}^{}\le C\,,\ \ \
\|\cc_\EPS\|_{L^\infty(I;L^2(\varOmega;\R^n))}^{}\le C\,,
\\
\label{est-c-eps.b}
&\Big\|\sum_{i=1}^nc_{i,\EPS}-1\Big\|_{L^\infty(I;L^2(\varOmega))}^{}\le C\sqrt\EPS\,,
\ \text{ and}
\\
\label{est-c-eps.c}
&\|\min(0,c_{i,\EPS})\|_{L^\infty(I;L^2(\varOmega))}^{}\le C\sqrt\EPS\ \ \
\text{ for  any }\ i=1,...,n\,,
\end{align}
 where the constants $C$ are independent of $\EPS$.   Note in
 addition that the projection $\bbP$ commutes with differentiation,
 namely $\nabla \bbP \bm\mu_\EPS = \bbP\nabla \bm\mu_\EPS$, where the
 latter projection is intended  column-wise. This fact and
 \eqref{est-mu-c-eps} in particular ensure  that
\begin{align}
 \|\bbP\bm\mu_\EPS\|_{L^2(I;H^1(\varOmega;\R^n))}^{} \le C\,,  \label{est-mu-c-eps2}
\end{align}
Moreover, using  \eqref{nabla-mu} written in the  equivalent  form
\begin{align}
  \nonumber
  \nabla\cc_\EPS&=\Big(\frac{[\phi_{\rm ref}^{\bm\xi} ]_{\cc\cc}'' 
    (J_\EPS,\cc_\EPS)\!}{J}+  \mathcal{P}_\EPS''(\cc_\EPS) \Big)^{-1} 
  \bigg(\nabla\bm\mu_\EPS 
  \\&\nonumber\ \ \ \ \ \ \ \ \ \
  -\frac{\!J_\EPS[\phi_{\rm ref}^{\bm\xi}]_{J\cc}''(J_\EPS,\cc_\EPS)
    {-}[\phi_{\rm ref}^{\bm\xi}]_{\cc}'(J_\EPS,\cc_\EPS)}{J^2}\nabla J_\EPS
  -\frac{[[\phi_{\rm ref}]_{\XX\cc}'']^{\bm\xi}(J_\EPS,\cc_\EPS)}{J_\EPS}\nabla\bm\xi\bigg)\,,
  \intertext{we obtain  the  estimate}
  &\|\nabla\cc_\EPS\|_{L^2(I{\times}\varOmega;\R^{n\times3})}^{}\le C\,
  \intertext{independently of $\EPS>0$. In addition,  \eqref{est-mu-c-eps} allow for using \eqref{Euler3-diff}
    written for $\EPS$-solution to obtain the estimate}
  &\Big\|\pdt{\cc_\EPS}\Big\|_{L^2(I;H^1(\varOmega;\R^n)^*)}^{}\le C\,.
\end{align}
\end{subequations}
Here, the bound on $\nabla \bm\mu_\EPS$ from \eqref{est-mu-c-eps} is used. 

\medskip\noindent{\it Step 6: convergence for $\EPS\to0$}.  
As in Step~4, we select a weakly* convergent subsequence and  find
$(\varrho,\vv,\bm\xi,\TTT J,\EEE\GRAVPOT,\cc,\bm\mu)$ such that
\begin{subequations}\label{eq:Regul.Sol}
\begin{align}
&\varrho_\EPS\to\varrho&&\text{weakly* in }\ L^\infty(I;W^{1,r}(\varOmega))
\,\cap\,W^{1,\TTT q\EEE}(I;L^r(\varOmega))\,,
\\&\vv_\EPS\to\vv&&\text{weakly* in }\ L^\infty(I;\TTT L^2\EEE(\varOmega;\R^3))
\,\cap\, 
\TTT L^q(I;W^{2,q}\EEE(\varOmega;\R^3))\,,
\\&\label{conv-of-xi}\bm\xi_\EPS\to\bm\xi&&\text{weakly* in }\
L^\infty(I;W^{2,r}(\varOmega;\R^3))\,\cap\,
W^{1,\TTT q\EEE}(I;L^{\TTT r\EEE}(\varOmega;\R^3))\,,
\\&
\TTT J_\EPS\to J\EEE&&\TTT\text{weakly* in }\
L^\infty(I;W^{1,r}(\varOmega))\,\cap\,W^{1,q}(I;L^{\TTT r\EEE}(\varOmega))\,,\EEE
\\&\GRAVPOT_\EPS\to\GRAVPOT&&\text{weakly* in }\
L^\infty(I;W^{1,r}(\varOmega))\,\cap\,W^{1,q}(I;L^{r}(\varOmega))\,,\EEE
{BV}(I;\TTT H^1\EEE(U))\,,
\\&\label{conv-of-c}
\cc_\EPS\to\cc&&\text{weakly in }\ L^2(I;H^1(\varOmega;\R^n))\,\cap\,
H^1(I;H^1(\varOmega;\R^n)^*)\,,\\
\label{conv-of-mu}
  &  \bbQ \bm\mu_\EPS \to \bm\mu &&\text{weakly in }\ 
 L^2(I;H^1(\varOmega;\R^n)) \, ,  
\end{align}
\end{subequations}
where $\bm\mu(t,\cdot)=\bbQ \bm\mu (t,\cdot)$ and the projection $\bbQ: L^2(\varOmega;\R^n) \to L^2(\varOmega;\R^n) $ is
defined via 
\[
\big(\bbQ\bm\mu\big)(x) = \bm\mu (x) - \left(\frac1{n|\varOmega|}
\int_\varOmega\bm\mu(y){\bm\cdot} \bm1\; \d y \right)\bm1\,.
\]
To justify convergence \eqref{conv-of-mu} we observe that 
$\int_\varOmega\big( |\nabla\bm\mu|^2{+}|\bbP\bm\mu|^2\big) {\rm d}
 \bm x  =0$ implies 
first, by using $\nabla\bm\mu=0$ that $\bm\mu(x)=\bm\eta$ for a constant vector
$\bm\eta\in \R^n$. Second, using $\bbP\bm\eta=0$  we have $\bm\eta=\alpha \bm1$
for a constant $\alpha\in \R$. Hence, we have $\bbQ\bm\mu=0$. Thus, by
repeating the classical compactness argument for showing Poincar\'e's
inequality, we find a constant $c_\bbQ >0 $ such that 
\[
\forall\, \bm\mu\in H^1(\varOmega;\R^n) \text{ with }\bbQ\bm\mu=0: \quad 
\int_\varOmega \big( |\nabla\bm\mu|^2 +|\bbP\bm\mu|^2\big) \;\d 
\bm x   \geq c_\bbQ \|
\bm\mu\|_{H^1(\varOmega;\R^n)}^2. 
\]
Thus, with \eqref{est-mu-c-eps} we obtain the a priori estimate 
\[
\| \bbQ \bm\mu_\EPS \|_{L^2(I;H^1(\varOmega;R^n))} \leq C 
\]
and \eqref{conv-of-mu} follows. 
 The Aubin-Lions lemma and  \eqref{conv-of-xi} and
\eqref{conv-of-c} yield the strong convergences
\begin{subequations}
  \begin{align}
    \bm \xi_\EPS \to \bm\xi\quad&\text{strongly in} \ \ C(I
                                  {\times} \overline \varOmega;\R^3)\,,\label{str_xi}\\
    J_\EPS \to J \quad&\text{strongly in} \ \ C(I
                        {\times} \overline \varOmega )\,,\label{str_J}\\
    \cc_\EPS \to \cc \quad&\text{strongly in} \ \
                            L^s(I{\times}\varOmega;\R^n) \ \ \forall 1\le
                            s<6\,. \label{str_c}
  \end{align}\label{str}
\end{subequations}

In addition to the  argument used for proving 
Theorem~\ref{thm-1}, we 
need  to pass to the limit in 
the semi-linear 
transport-and-diffusion equation \eqref{Euler3-diff} formulated
weakly.  This is however straightforward, so that we omit details. 

The  only  remaining  issue is to pass to the limit  in
\eqref{Euler3-diff++} written
 at level  $\EPS$  
in the form of the variational inequality  
\begin{align}
&\int_0^T\!\!\!\int_\varOmega\mathcal{P}_\EPS( \cc_\EPS
                 )+
\Big(\bm\mu_\EPS
-\frac{[\phi_{\rm ref}^{\bm\xi_\EPS}]_\cc'(J_\EPS,\cc_\EPS)}{J_\EPS}\Big)
{\cdot}( \widetilde \cc_\EPS  - \cc_\EPS ) \,\d\xx \;\!\d t
 \leq  \int_0^T\!\!\!\int_\varOmega\mathcal{P}_\EPS( \widetilde \cc_\EPS)\,\d\xx \;\!\d t
\label{w-sln-mu++}
\end{align}
for $i=1,...,n$, for  all  $ \wt\cc_\EPS  \in
L^2(I{\times}\varOmega;\R^n)$.  
 The bounds (\ref{est-c-eps}b,c) entail that the limit $\cc$ is
valued in $\triangle_1^+$ almost everywhere in $I \times \varOmega$. 
We aim at proving
that $\cc$ 
satisfies 
\begin{align}
\int_0^T\!\!\!\int_\varOmega
\Big(\bm\mu 
  - 
  \frac{[\phi_{\rm ref}^{\bm\xi}]_\cc'(J,\cc)}{J}\Big)
{\cdot}(\widehat\cc-\cc)\,\d\xx \;\!\d t \leq  0
\label{w-sln-mu+++}\end{align}
for  all  $\widehat\cc\in L^2(I{\times}\varOmega;\R^n)$ valued in
$\triangle_1^+$ a.e.\ on $I{\times}\varOmega$.
 We first observe that \eqref{est-c-eps.b} allows to
write 
\[
\bm c_\EPS(t,x)= \displaystyle\frac{1+\alpha_\EPS (t)}{n} \bm1 + (\bbQ \bm c_\EPS\big) (t,x), \quad
\text{where } \|\alpha_\EPS \|_{L^\infty(I)} \leq C\sqrt{\EPS}.
\]
To construct good test function $\wt\cc_\EPS$ for \eqref{w-sln-mu++} 
we fix a small $\delta>0$ and choose $\widehat\cc$ with
\[
\widehat\cc(t,x) \in \triangle^+_1\text{ and }\ \min \widehat c_i(t,x)\geq
\delta \quad \text{and set } \ \wt\cc_\EPS= \displaystyle\frac{
  \alpha_\EPS (t)}{n} \bm1 +
\widehat\cc. 
\] 
In particular, we have $\wt\cc_\EPS- \cc_\EPS= \bbQ\big(\widehat\cc -
{\bm c}_\EPS  \big) $.  By using the weak convergences
(\ref{eq:Regul.Sol}e,f)
and the strong convergences
\eqref{str},  we obtain 
\begin{align*}
&\int_0^T\!\!\!\int_\varOmega
\Big(\bm\mu_\EPS -\frac{[\phi_{\rm ref}^{\bm\xi_\EPS}]_\cc'(
                 J_\EPS,\cc_\EPS)}{ J_\EPS}\Big)
{\cdot}(\wt\cc_\EPS-\cc_\EPS)\,\d\xx \;\!\d t 
=
\int_0^T\!\!\!\int_\varOmega
\Big(\bbQ\bm\mu_\EPS-\bbQ\frac{[\phi_{\rm ref}^{\bm\xi_\EPS}]_\cc'(
                 J_\EPS,\cc_\EPS)}{ J_\EPS}\Big)
{\cdot}(\widehat\cc-\cc_\EPS)\,\d\xx \;\!\d t 
\\
&\quad \longrightarrow 
\int_0^T\!\!\!\int_\varOmega
\Big(\bm\mu-  \bbQ  \frac{[\phi_{\rm ref}^{\bm\xi}]_\cc'(J,\cc)}{J}\Big)
{\cdot}(\widehat\cc-\cc)\,\d\xx \;\!\d t   = \int_0^T\!\!\!\int_\varOmega
\Big(\bm\mu-   \frac{[\phi_{\rm ref}^{\bm\xi}]_\cc'(J,\cc)}{J}\Big)
     {\cdot}\bbQ(\widehat\cc-\cc)\,\d\xx \;\!\d t   \\
  &\quad  = \int_0^T\!\!\!\int_\varOmega
\Big(\bm\mu-   \frac{[\phi_{\rm ref}^{\bm\xi}]_\cc'(J,\cc)}{J}\Big)
{\cdot}(\widehat\cc-\cc)\,\d\xx \;\!\d t\,,
\end{align*}
where we have also used  $\bbQ(\widehat\cc-\cc)= \widehat\cc-\cc $,
coming from the fact that $\cc, \, \widehat \cc \in
\triangle_1^+$.  

Moreover, our construction guarantees $\widetilde
c_{\EPS,i}(t,x)=\alpha_\EPS(t) /n 
  + \widehat c_i(t,x)\geq -C\sqrt\EPS + \delta\geq 0$ for sufficiently small
  $\EPS>0$. Hence, we have $\mathcal P_\EPS(\wt\cc_\EPS)=
  \big(\alpha_\EPS{+}\widehat\cc{\bm\cdot}\bm1-1\big)^2/(2\EPS)=
  \alpha_\EPS^2/(2\EPS) $ since $\widehat\cc \in
  \triangle^+_1$. With this, $\mathcal P_\EPS(\cc_\EPS)\geq
  (\cc_\EPS{\cdot}\bm1 -1)^2/(2\EPS)$, and $\alpha_\EPS (t)=\frac1{|\varOmega|}
  \int_\varOmega \big(\cc_\EPS(t,y){\cdot} \bm 1-1\big)\,\d y$ we find 
\begin{align*}
&\int_0^T\!\!\!\int_\varOmega \big(\mathcal P_\EPS(\cc_\EPS)
- \mathcal P_\EPS(\wt\cc_\EPS)\big) \:\d x \;\! \d t \geq 
\frac1{2\EPS} \int_0^T\!\!\!\int_\varOmega \Big(
  \big(\cc_\EPS{\cdot}\bm1 -1\big)^2 - 
 \big(\alpha_\EPS\big)^2 \Big) \:\d x  \;\!\d t\ \geq  \ \frac1{2\EPS} \int_0^T
 0 \,\d t =0.
\end{align*}

Hence,  by collecting all terms on the left-hand side, we can pass to the
$\liminf$ $\EPS\to 0$ in \eqref{w-sln-mu++} and obtain \eqref{w-sln-mu+++} 
for all $\widehat \cc$ satisfying $\min \widehat
c_i(t,x)\geq \delta$. Since $\delta>0$  was arbitrary, the desired variational
inequality \eqref{w-sln-mu+++} holds for all test functions.

\medskip\noindent{\it Step 7: energy-dissipation balance}.
In addition to the argumentation used in the proof of  Theorem~\ref{th:WeakSol}(ii),
 we now use that $\pdt{}\cc+(\vv{\cdot}\nabla)\cc
\in L^2(I;H^1(\varOmega;\R^n)^*)$.  Note also that 
$(\vv{\cdot}\nabla)\cc$ lies in $L^2(I;L^{2}(\varOmega;\R^n))$
and can be tested by $\bm\mu$ and integrated by parts. Here the indeterminacy
of $\bm\mu$ with  respect  to spatially constant multiples of $\bm1$ does not
matter, because from $\cc(t,x) \in \triangle^+_1$ we have $\cc(t,\xx) \cdot
\bm1 =0$ in $I{\times} \varOmega$.  Thus, the energy balance follows as for
Theorem~\ref{th:WeakSol}(ii). 
\end{proof}

\begin{remark}[{\sl Composition-dependent mass density}]\label{rem=mass}\upshape
 It would be desirable to make the  mass density depends also on composition, cf.\
e.g.\ \cite[Ch.2]{Gery19INGM}. In our present modeling level, this would mean
$\rho_{\rm ref}=\rho_{\rm ref}(\XX,\cc)$ and  hence,  by  
\eqref{rho=rho0/detF}  we would have  $\varrho =\rho_{\rm
  ref}^{\bm\xi}(\cc)/J$.
 Repeating the calculation in \eqref{DT-det} and \eqref{transport-xi}  
continuity equation \eqref{cont-eq+}  would be extended to 
\begin{align}
\DT\varrho+({\rm div}\,\vv)\varrho=\frac{[\rhoREF^{\bm\xi}]_\cc'(\cc)}{J}{\cdot}\DT\cc\,.
\label{towards-cont-eq+}
\end{align}
 Thus, the principle of mass conservation would be violated, and further
mathematical difficulties would arise in the  kinetic energy and the
gravitational energy \eqref{calculus-selfgravit}.  Hence, a proper modeling of
a concentration-dependent mass density would need a truly multiphase modeling
that exceeds beyond the Eckart-Prigogine approximation and is left to future
research.  
\end{remark}

\begin{remark}[{\sl An alternative approach}]\label{rem-M}\upshape
Keeping  the  sum of diffusion fluxes  to $0$ 
can  also be 
achieved
without  directly constraining $\sum_{j=1}^n c_j= 1$  
by tuning the mobility matrix
to satisfy $\sum_{j=1}^n \bbM_{ij}(\XX,J,\cc)=0$ ,  which is the ``local
model'' in the sense of \cite{OttE97TDIF}.  In this case,
condition $\sum^n_{i=1}c_i=1$ is kept during the evolution if it holds  at
the initial time.  
This can be seen by summing up \eqref{Euler3-diff} for $i=1,...,n$, which gives
$\pdt{}(\sum^n_{i=1}c_i)={\rm div}\big(\sum_{i,j=1}^n
\bbM_{ij}^{\bm\xi}(J,\cc)\nabla\mu_j \big) +\sum_1^n r_i^{\bm\xi}(J,\cc))={\rm
  div}\,\bm0+0=0$. This allows to us avoid the constraint $\sum^n_{i=1}c_i=1$
from $\bm N(\cdot)$ in \eqref{Euler3-diff+}.  For equal mobilities of each
chemical components,  possibly  dependent on pressure and local
composition, denoted by $m=m(\XX,J,\cc)$, one usually considers the so-called
Maxwell-Stefan mobility matrix
\begin{align}\label{M}
\bbM(\XX,J,\cc)=m(\XX,J,\cc)\big({\rm diag}(\cc)-\cc{\otimes}\cc\big)\,,
\end{align}
cf.\ e.g.\ \cite{Giov99MFM, BotDru23SCTD}. This matrix  has the kernel
$\bm1$ and hence is only positive semidefinite; moreover it degenerates further
for $c_i\approx 0$.  This makes its usage analytically  more difficult
than our simplified model with a  general symmetric positive matrix as
imposed in \eqref{Euler-fluid-selfgravit+}.
\end{remark}

\begin{remark}[\emph{General reaction stoichiometry}]
\label{rem:Stoich}\upshape
For general reaction systems the assumption that $\bbK$  in
\eqref{ass-r}  is positive definite on
the orthogonal complement of $\bm1$ may be too restrictive. In general, one has
a stoichiometric subspace $\bm S\subset \R^n$ such that $\rr \in \bm
S$. Defining $\bbP_S:\R^n\to \bm S\subset \R^n$ to be the orthogonal projection
and $\bbQ_S=I-\bbP_S$ the complementing projection, one can then assume that
there is a symmetric reaction matrix $\bbK$ such that
$\mu\cdot \bbK(\XX,J,\cc) \mu \geq \epsilon |\bbP_S \mu|^2$. The above analysis
can easily be generalized to this case, if we use that
$\bbQ_S \int_\varOmega \cc(t,\xx)\,\d\xx$ is conserved along solutions, that
$\bbP_S\bm\mu$ is controlled by the dissipation, and that we can assume
$\bbQ_S\int_\varOmega \bm\mu(t,\xx)\,\d\xx=0$ without loss of generality.  
\end{remark}

\bigskip

{\renewcommand{\baselinestretch}{0.9}
  \small 

\noindent{\it Acknowledgments.}
 A.M.\ was partially supported by Deutsche Forschungsgemeinschaft (DFG) through
the Berlin Mathematics Research Center MATH+ (EXC-2046/1, project ID
390685689, subproject \emph{DistFell}).
U.S.\ received partial
support  from the Austrian Science Fund (FWF) projects F\,65, I\,4354,
I\,5149, P\,32788. 
 T.R.\ was partially supported by  the Czech Sci.\ Foundation (CSF/DFG
project GA22-00863K) with 
the institutional support RVO:61388998 (\v CR).

\baselineskip=10pt

\newcommand{\MR}[1]{}

}


\newcommand{\MR}[1]{}
\newcommand{\etalchar}[1]{$^{#1}$}
\providecommand{\bysame}{\leavevmode\hbox to3em{\hrulefill}\thinspace}

\end{document}